\documentclass{amsart}

\usepackage[T1]{fontenc}
\usepackage[latin1]{inputenc}
\usepackage[english]{babel}
\usepackage{amsmath,amssymb}
\IfFileExists{mathpazo.sty}{
  \usepackage{mathpazo}
}{}
\usepackage{mathrsfs}
\usepackage{andsupp}
\usepackage{amscd}
\usepackage{boxedminipage}
\usepackage{ifpdf}
\ifpdf
\usepackage[line,cmtip,all,matrix]{xy}
\else
\usepackage[line,cmtip,all,matrix,dvips]{xy}
\fi

\renewcommand{\epsilon}{\varepsilon}
\renewcommand{\rho}{\varrho}

\newcommand{\Sp}{\operatorname{S}}

\newcommand{\h}{\mathscr{H}}

\newtheorem{proposition}{Proposition}[section]
\newtheorem{lemma}{Lemma}[section]
\newtheorem{theorem}{Theorem}[section]
\newtheorem{corollary}{Corollary}[section]
\newtheorem{definition}{Definition}[section]
\newtheorem{claim}{Claim}[subsection]

\newtheorem{remark}{Remark}[section]

\newcommand{\e}[1]{\mathbf #1}
\renewcommand{\varinjlim}{\mathop{\fam0 \stackrel{\textstyle Lim}{\textstyle \longrightarrow}}}

\newcommand{\Dim}{\mathop{\fam0 Dim}\nolimits}

\newcommand{\Hom}{{\hbox{\rm Hom}}}

\newcommand{\Bosonnormalord}{\stackrel{\textstyle \circ}{ \circ}}

\newcommand{\tr}{\mathop{\fam0 Tr}\nolimits}

\newcommand{\ad}{\mathop{\fam0 ad}\nolimits}

\newcommand{\id}{\mathop{\fam0 Id}\nolimits}
\newcommand{\Lie}[1]{\mbox{\sf #1}}

\newcommand{\Id}{\mathop{\fam0 Id}\nolimits}

\newcommand{\R}{\mathop{\fam0 {\mathbb R}}\nolimits}

\newcommand{\bC}{\mathop{\fam0 {\mathbb C}}\nolimits}

\newcommand{\Z}{{\mathbb Z}}

\newcommand{\Tr}{{\rm tr}}

\newcommand{\ra}{\mathop{\fam0 \rightarrow}\nolimits}

\def\l{\lambda}\def\m{\mu}
\def\L{\Lambda}

\newcommand{\fg}{\mathop{\fam0 {\mathfrak g}}\nolimits}

\newcommand{\Si}{\Sigma}
\newcommand{\Sib}{{\mathbf \Sigma}}

\def\Vdag{{\mathcal V}^{\dagger}}

\def\ad{\operatorname{ad}}
\def\tr{\operatorname{tr}}

\def\:{:}

\def\Hom{\operatorname{Hom}}

\def\Bosonnormalordconstruction#1{\vcenter{\hbox{\ooalign{%
\raise.8ex\hbox{$#1\circ$}\crcr\lower.8ex\hbox{$#1\circ$}}}}}

\newcommand{\morf}[4][\longrightarrow]{#2 \colon #3#1#4}

\newcommand{\abs}[1]{\lvert#1\rvert}
\newcommand{\bo}[1][\lambda]{{\Box_{#1}}}

\newcommand{\hsu}{H^{\operatorname{SU}(N)}_K}
\newcommand{\inner}[1]{\langle#1\rangle}
\newcommand{\V}{\mathcal{V}}

\DeclareMathOperator{\End}{End}

\def\Bosonnormalord{\,\lower.8ex \hbox{$\circ$} \llap{\raise.8ex\hbox{$\circ$}} \,}
\def\normalord{\,\lower.8ex \hbox{$\cdot$} \llap{\raise.8ex\hbox{$\cdot$}} \,}
\def\Hom{\mathop{\rm Hom}\nolimits}

\def\ra{\rightarrow}

\def\ad{\mathop{\rm ad}\nolimits}

\def\tr{\mathop{\rm Tr}\nolimits}

\def\bC{{\mathbb C}}

\def\bN{{\mathbb N}}

\def\N{{\mathcal N}}

\def\gg{\mathfrak{g}}

\def\gh{\mathfrak{h}}

\def\Id{\mathop{\rm Id}\nolimits}

\def\H{{\mathcal H}}

\begin{document}

\title[The WRT-TQFT from Conformal Field Theory]{Construction of the
Witten-Reshetikhin-Turaev TQFT from conformal field theory}

\author{J{\o}rgen Ellegaard Andersen}
\address{Center for Quantum Geometry of Moduli Spaces\\
        University of Aarhus\\
        DK-8000, Denmark}
\email{andersen{\@@}imf.au.dk}

\author{Kenji Ueno}
\address{Department of Industrial and System Engineering\\ Faculty of Science and Engineering\\
Hosei University\\ Koganei, Tokyo\\ 184-0002 Japan}
\email{ueno{\@@}math.kyoto-u.ac.jp}

\begin{abstract}  In \cite{AU2} we constructed the {\em vacua modular
functor} based on the sheaf of vacua theory developed in \cite{TUY}
and the abelian analog in \cite{AU1}. We here provide an explicit
isomorphism from the modular functor underlying the skein-theoretic
model for the Witten-Reshetikhin-Turaev TQFT due to Blanchet, Habbeger,
Masbaum and Vogel to the vacua modular functor. This thus provides a
geometric construction of the TQFT first proposed by Witten and constructed first by
Reshetikhin-Turaev from the quantum group $U_q(\Lie{sl}(N))$.
\end{abstract}

\thanks{Supported in part by the center of excellence grant ``Center
  for quantum geometry of Moduli Spaces'' from the Danish National
  Research Foundation and by Grant-in-Aid for Scientific Research
  No.19340007 from JSPS (Japan Society for the Promotion of Science)}

\maketitle



\section{Introduction}

This is the main paper in a series of four papers (the previous three being \cite{AU1}, \cite{AU2}, \cite{AU3}), where we
provide a geometric construction of modular functors and
topological quantum field theories (TQFT's) from conformal field theory
building on the constructions of Tsuchiya, Ueno and Yamada in \cite{TUY} and \cite{Ue2} and
Kawamoto, Namikawa, Tsuchiya and Yamada in \cite{KNTY}. In this paper we provide an explicit isomorphism of the modular
functors underlying the Witten-Reshetikhin-Turaev TQFT for
$U_q(\Lie{sl}(N))$ \cite{RT1,RT2,TW2} and the {\em vacua modular functors} constructed in \cite{AU2}, based on the conformal field theory
for the Lie algebra $\Lie{sl}(N)$ constructed in \cite{TUY} and \cite{Ue2} and the one dimensional ghost theory constructed in
\cite{KNTY}.  We use the skein theory
approach to the Witten-Reshetikhin-Turaev TQFT of Blanchet, Habegger,
Masbaum and Vogel \cite{BHMV1}, \cite{BHMV2} and \cite{Bl1} together with work of Wenzl \cite{Wzl1} and Kanie \cite{Kanie1} to set up this
isomorphism. Since the modular functor determines the TQFT
uniquely, this therefore also provides a geometric construction of
the Witten-Reshetikhin-Turaev TQFT. That there should be such an isomorphism is a well established conjecture which is due to Witten, Atiyah and Segal (see e.g. \cite{W}, \cite{Atb} and \cite{Se}).

Let us now outline our construction of the above mentioned isomorphism between the two theories.
In \cite{AU1} we described how to reconstruct the rank one ghost theory first introduced by Kawamoto, Namikawa, Tsuchiya and Yamada in \cite{KNTY}
from the the point of view of \cite{TUY} and \cite{Ue2}.
In \cite{AU2} we described how one combines the work of Tsuchiya, Ueno and Yamada (\cite{TUY} and \cite{Ue2}) with \cite{AU1} to construct the vacua modular functor for each simple Lie algebra and a positive integer $K$ called the level. See also \cite{Ue4}. Let us here denote the theory we constructed for  the Lie algebra $\Lie{sl}(N)$ and level $K$ by $\Vdag_{N,K}$. 
We recall that a modular functors is a functor from a certain category of extended labeled marked surfaces (see section \ref{AxiomsMF}) to the category of finite dimensional vector spaces. The functor is required to satisfy Walker's topological version \cite{Walker} of Segal's axioms for a modular functor \cite{Se} (see section \ref{AxiomsMF}). Note that we do not consider the duality axiom as part of the definition of a modular functor. We consider the duality axioms as extra data. For modular functors which satisfies the duality axiom, we say that it is a modular functor with duality.
In \cite{Bl1} Blanchet constructed a modular tensor category which we will here denote $\hsu$ (see section \ref{cha:notes-part-one}). It is constructed using skein theory and one can build a modular functor and a TQFT from this category following either the method of \cite{BHMV2} or \cite{Tu}. We denote the resulting modular functor $\V_K^{SU(N)}$.
It is easy to check that the two modular functors $\Vdag_{N,K}$ and $\V_K^{SU(N)}$ have the same label set $\Gamma_{N,K}$.
In this paper we explicitly construct an isomorphism between these to modular functors.

\begin{theorem}\label{Main}
There is an isomorphism of modular functors
\[I_{N,K} : \V_K^{SU(N)} \ra \Vdag_{N,K},\]
i.e. for each extended labeled marked surface $(\Sib, \lambda)$ we have an isomorphism of complex vector spaces
\[I_{N,K}(\Sib,\lambda) : \V_K^{SU(N)}(\Sib,\lambda) \ra \Vdag_{N,K}(\Sib,\lambda),\]
which is compatible with all the structures of a modular functor.
\end{theorem}

The choices involved in our construction of this isomorphism is only a non-zero complex scalar  for each element in the label set $\Gamma_{N,K}$ (see Definition \ref{swfix}). If we multiply the choice for $I_{N,K}$ by the non-zero complex scalars $c_\lambda$ for each label $\lambda\in \Gamma_{N,K}$, then the new isomorphism $\tilde I_{N,K}$ for any labeled marked surface $(\Sib,\lambda)$ is given by
$$\tilde I_{N,K}(\Sib,\lambda)  = \prod_{p\in P} c_{\lambda(p)} I_{N,K}(\Sib,\lambda),$$
where $P$ is the set of labelled marked points of $\Sib$. Hence up to these simple automorphisms, which any modular functor has, the isomorphism is independent of these choices.

The main idea behind the construction of $I_{N,K}$ is to use the GNS construction applied to the infinite Hecke algebra with respect to the relevant Markov traces as was first done by Jones \cite{J-4} and Wenzl \cite{Wzl1}. On the skein theory side, we identify the usual purification construction in terms of the GNS construction. This allows us to show that the resulting representations of the Hecke algebras are isomorphic to Wenzl's representations. Analyzing the properties of this isomorphism further we find that it determines isomorphism from certain morphism spaces in the category $\hsu$ to 
Wenzl's representation spaces. On the vacua side, we know by the results of Kanie \cite{Kanie1} (see also \cite{Ue3}) that the space of vacua gives a geometric construction of Wenzl's representations.  In fact Kanie constructs explicit isomorphisms between certain space of vacua and Wenzl's representations. Combining the isomorphisms arising from the above mentioned analysis of the GNS-constructions with Kanie's isomorphisms, we arrive at the very important isomorphisms between certain morphism spaces of the category $\hsu$ and certain spaces of vacua given in Definition \ref{hsuglabel}. From the very way these isomorphisms are constructed, they provide an identification of certain factorisations on the skein theory side to those on the vacua modular functor side as stated in Theorem \ref{decompI}. By following Turaev's construction of a modular functor from a modular tensor category, we build the the modular functor $ \V_K^{SU(N)} $ from the modular category $\hsu$ and the isomorphism from Definition \ref{hsuglabel} now determines all the needed isomorphism between the vector spaces $ \V_K^{SU(N)} $ and $\Vdag_{N,K}$ associates to all labeled marked surfaces. Theorem \ref{decompI} is key in showing that this is indeed an isomorphism of modular functors in genus zero. The main result of \cite{AU3} then implies that the isomorphisms of modular functors (without duality) in genus zero can be extended to surfaces of all genus. It is only after Theorem \ref{Main} has been established that we get the structure of duality and unitarity on the vacua modular functor $ \Vdag_{N,K}$, since we can use our isomorphism $I_{N,K}$ to transfer these two structures from $\V_K^{SU(N)}$ to $ \Vdag_{N,K}$.  In fact at present we do not have a geometric construction of a duality structure for $ \Vdag_{N,K}$. 

We have the following geometric application of our construction.

\begin{theorem}\label{bofvunitary}
The connections constructed in the bundle of vacua for any holomorphic family of labeled marked curves given in \cite{TUY} preserves projectively a unitary structure which is projectively compatible with morphism of such families.
\end{theorem}
This Theorem is an immediate corollary of our main Theorem \ref{Main}. By definition $\Vdag_{N,K}(\Sib,\lambda)$ is the covariant constant sections 
of the bundle of vacua twisted by a fractional power of a certain ghost theory over Teichm\"{u}ller space as described in \cite{AU2}. Using the isomorphism $I_{N,K}$ from our main Theorem \ref{Main}, we transfer the unitary structure on $\V_K^{SU(N)}(\Sib,\lambda)$ to the bundle of vacua over Teichm\"{u}ller space. Here we have used the preferred section of the ghost theory, to transfer the unitary structure to the bundle of vacua (see \cite{AU2}). Since the unitary structure on $\V_K^{SU(N)}(\Sib,\lambda)$ is invariant under the extended mapping class group, the induced unitary structure on the bundle of vacua will be projectively invariant under the action of the mapping class group. But since the bundle of vacua for any holomorphic family naturally is isomorphism to the pull back of the bundle of vacua over Teichm\"{u}ller space, we get the stated theorem. As a further application we get that

\begin{theorem}\label{Hitchinunitary}
The Hitchin connections constructed in the bundle $\H^{(K)}$ over Teichm\"{u}ller space, whose fiber over an algebraic curve, representing a point in Teichm\"{u}ller space, is the geometric quantization at level $K$ of the moduli space of semi-stable bundles of rank $N$ and trivial determinant over the curve, projectively preserves a unitary structure which is projectively preserved by the mapping class group.
\end{theorem}
This is an immediate corollary of Theorem \ref{bofvunitary} and then the theorem by Laszlo in \cite{La1}, which provides a projective isomorphism of the bundle $\H^{(K)}$ with its Hitchin connection \cite{Hi} and then the bundle of vacua with the TUY-connections over Teichm\"{u}ller space. In the paper \cite{AUnitary} this unitary structure is constructed explicitly for $N=2$.  We also get the following corollary

\begin{corollary}
The projective monodromy of the Hitchin connection contains elements of infinite order for $N$ and $k\notin \{1,2,4,8\}$.
\end{corollary}

Masbaum proved the corresponding result for the Witten-Reshetikhin-Turaev Theory for $N=2$ in \cite{MFI}. This theorem is therefore an immediate corollary of Theorem \ref{Main} and the results of \cite{AMU}, which shows that Masbaum's arguments for $N=2$ implies the statement for any $N$ for the Witten-Reshetikhin-Turaev Theory. Prior to this, Funar proved, in the case of $N=2$, that the monodromy group was infinite for the same levels for the Witten-Reshetikhin-Turaev theory in \cite{Fu}. For a purely algebraic geometric proof of this result for $N=2$, see \cite{LPS}.

A further corollary is that the projective space of projective flat sections of the bundle $\H^{(K)}$ over Teichm\"{u}ller space of any closed surface $\Sigma$ is isomorphic to the space ${\mathbb P}\V_K^{SU(N)}(\Sigma)$. Therefore we are allowed to use the geometric quantization of the moduli space of flat $SU(N)$-connections to study the Witten-Reshetikhin-Turaev TQFT's. We have already provided a number of such applications, see e.g. \cite{A5, AMas, A1, AG, A2, A3, A4, ALG,  AB, AH, APropT, Agf}.

The paper is organized as follows. 
In section \ref{AxiomsMF} we briefly recall the theory of modular functors. 
Section \ref{cha:notes-part-one} briefly describes, following Blanchet \cite{Bl1}, the modular tensor category $\hsu$ using skein theory. In the following section \ref{RTmf}, we define the modular functor 
$\V_K^{SU(N)}$ which gives the skein-theoretic model for Witten-Reshetikhin-Turaev TQFT \cite{RT1,RT2,TW2}
due to Blanchet, Habbeger, Masbaum and Vogel \cite{BHMV1}, \cite{BHMV2}. 
In section \ref{sec:rep-Hecke} we recall Jones and Wenzl work on the representation theory of the Hecke algebra and establish the needed relation between the skein theory representations and Wenzl's representations of the Hecke algebra.
The genus zero part of the isomorphism $I_{N,K}$ between the 
modular functors  $\V_K^{SU(N)}$ and $\Vdag_{N,K}$ is provided in section \ref{cha:Identification}. It is recalled in section \ref{SMatrixhg} how we in the paper \cite{AU3} proved that such a genus zero isomorphism automatically extends to a full isomorphism of modular functors, thus establishing our main Theorem \ref{Main}.

We thank Christian Blanchet, Yukihiro Kanie, Gregor Masbaum, Akihiro Tsuchiya and Yasuhiko Yamada for valuable discussions and the referee for several helpful comments which improved the paper.

\section{Modular functors}\label{AxiomsMF}

\subsection{The axioms for a modular functor}

We shall in this section give the axioms for a modular functor.
These are due to G. Segal and appeared first in \cite{Se}. We
present them here in a topological form, which is due to K. Walker
\cite{Walker}. See also \cite{Kont1, Grove}. We note that similar, but
different, axioms for a modular functor are given in \cite{Tu} and
in \cite{BB}. The authors are not aware of a proof of the equivalence of these definitions of
a modular functor. However, we will not need it in this paper.

For a closed oriented surface $\Si$ of genus $g$ we have the
non-degenerate skew-symmetric intersection pairing
\[(\cdot,\cdot) : H_1(\Si,\Z) \times H_1(\Si,\Z) \ra \Z.\]
Suppose $\Si$ is connected. In this case a Lagrangian subspace $L\subset H_1(\Si,\Z)$
is by definition a subspace, which is maximally isotropic with respect to the
intersection pairing. If $\Si$ is not connected, then $H_1(\Si,\Z) =
\oplus_i H_1(\Si_i,\Z)$, where $\Si_i$ are the connected
components of $\Si$. By definition a Lagrangian subspace is in this
paper a subspace of the form $L = \oplus_i L_i$, where each $L_i\subset
H_1(\Si_i,\Z)$ is Lagrangian. 

For any real vector space $V$, we define $PV = (V-\{0\})/\R{}_+.$

\begin{definition} \label{msurface}

A {\em marked surface\/} $ {\Sib} = (\Si, P, V, L)$ is an oriented
closed smooth surface $\Si$ with a finite subset $P \subset \Si$
of points with projective tangent vectors $V\in \sqcup_{p \in
P}PT_{p}\Si$ and a Lagrangian subspace $L \subset H_1(\Si,\Z)$.
\end{definition}

\begin{definition} \label{mmorphism}

A {\em morphism\/} $\e f : {\Sib}_1 \to {\Sib}_2$ of marked
surfaces ${\Sib}_i = (\Si_i,P_i,V_i,L_i)$ is an isotopy class of
orientation preserving diffeomorphisms $f : \Si_1 \to \Si_2$ that
maps $(P_1,V_1)$ to $(P_2,V_2)$ together with an integer $s$.
Hence we write $\e f = (f,s)$.
\end{definition}

Let $\sigma$ be Wall's signature cocycle for triples of Lagrangian
subspaces of $H_1(\Si,\R{})$ (See \cite{Wall}).

\begin{definition} \label{composition}
Let $\e f_1 = (f_1,s_1) : {\Sib}_1 \to {\Sib}_2$ and $\e f_2 =
(f_2,s_2) : {\Sib}_2 \to {\Sib}_3$ be morphisms of marked surfaces
${\Sib}_i = (\Si_i,P_i,V_i,L_i)$ then the {\it composition\/} of
$\e f_1$ and $\e f_2$ is $$ \e f_2 \e f_1 = (f_2 f_1, s_2 + s_1 -
\sigma((f_2f_1)_*L_1, f_{2*}L_2,L_3)). $$
\end{definition}

With the objects being marked surfaces and the morphism and their
composition being defined as in the above definition, we have
constructed
 the category of marked surfaces.

The mapping class group $\Gamma({\Sib})$ of a marked surface
${\Sib} = (\Si,L)$ is the group of automorphisms of ${\Sib}$. One
can prove that $\Gamma({\Sib})$ is a central extension of the
mapping class group $\Gamma(\Si)$ of the surface $\Si$ defined by
the 2-cocycle $c : \Gamma({\Sib}) \to \mathbb Z$, $c(f_1,f_2) =
\sigma((f_1f_2)_*L,f_{1*}L,L)$. One can also prove that this
cocycle is equivalent to the cocycle obtained by considering
two-framings on mapping cylinders (see \cite{At1}).

\begin{definition} \label{disjunion}
The operation of {\em disjoint union of marked surfaces} is $$
(\Si_1,P_1,V_1,L_1)
 \sqcup (\Si_2,P_2,V_2,L_2) = (\Si_1 \sqcup \Si_2,P_1 \sqcup P_2,V_1\sqcup V_2,L_1 \oplus
L_2). $$

Morphisms on disjoint unions are accordingly $(f_1,s_1) \sqcup
(f_2,s_2) = (f_1 \sqcup f_2,s_1 + s_2)$.
\end{definition}

We see that disjoint union is an operation on the category of
marked surfaces.

\begin{definition}\label{or}
Let ${\Sib}$ be a marked surface. We denote by $- {\Sib}$ the
marked surface obtained from ${\Sib}$ by the {\em operation of
reversal of the orientation}. For a morphism $\e f = (f,s) :
{\Sib}_1 \to {\Sib}_2$ we let the orientation reversed morphism be
given by
 $- \e f = (f,-s) : -{\Sib}_1 \to -{\Sib}_2$.
\end{definition}

We also see that orientation reversal is an operation on the
category of marked surfaces. Let us now consider glueing of marked
surfaces.

Let $(\Si, \{p_-,p_+\}\sqcup P,\{v_-,v_+\}\sqcup V,L)$ be a marked
surface, where we have selected an ordered pair of marked points
with projective tangent vectors 
$$((p_-,v_-),(p_+,v_+)),$$ 
at which
we will perform the glueing.

Let $c : P(T_{p_-}\Si) \ra P(T_{p_+}\Si)$ be an orientation
reversing projective linear isomorphism such that $c(v_-) = v_+$.
Such a $c$ is called a {\em glueing map} for $\Si$. Let
$\tilde{\Si}$ be the oriented surface with boundary obtained from
$\Si$ by blowing up $p_-$ and $p_+$, i.e.
\[\tilde{\Si} = (\Si -\{p_-,p_+\})\sqcup P(T_{p_-}\Si)\sqcup P(T_{p_+}\Si),\]
with the natural smooth structure induced from $\Si$. Let now
$\Si_c$ be the closed oriented surface obtained from $\tilde{\Si}$
by using $c$ to glue the boundary components of $\tilde{\Si}$. We
call $\Si_c$ the glueing of $\Si$ at the ordered pair
$((p_-,v_-),(p_+,v_+))$ with respect to $c$. 

Let now $\Si'$ be the topological space obtained from $\Si$ by
identifying $p_-$ and $p_+$. We then have natural continuous maps
$q : \Si_c \ra \Si'$ and $n : \Si \ra \Si'$. On the first homology
group $n$ induces an injection and $q$ a surjection, so we can
define a Lagrangian subspace $L_c \subset H_1(\Si_c,\Z)$ by $L_c =
q_*^{-1}(n_*(L))$. We note that the image of $P(T_{p_-}\Si)$ (with
the orientation induced from $\tilde{\Si}$) induces naturally an
element in $H_1(\Si_c,\Z)$ and as such it is contained in $L_c$.

\begin{remark}{\em \label{remarkglue2}
If we have two glueing maps $c_i : P(T_{p_-}\Si) \ra
P(T_{p_+}\Si),$ $i=1,2,$ we note that there is a diffeomorphism
$f$ of $\Si$ inducing the identity on
$(p_-,v_-)\sqcup(p_+,v_+)\sqcup(P,V)$ which is isotopic to the
identity among such maps, such that $(df_{p_+})^{-1} c_2 df_{p_-}
= c_1$. In particular $f$ induces a diffeomorphism $f : \Si_{c_1}
\ra \Si_{c_2}$ compatible with $f : \Si \ra \Si$, which maps
$L_{c_1}$ to $L_{c_2}$. Any two such diffeomorphims of $\Si$
induces isotopic diffeomorphims from $\Si_1$ to
$\Si_2$.}\end{remark}

\begin{definition} \label{glueing}
Let ${\Sib} = (\Si, \{p_-,p_+\}\sqcup P,\{v_-,v_+\}\sqcup V,L)$ be
a marked surface. Let $$c : P(T_{p_-}\Si) \ra P(T_{p_+}\Si)$$ be a
glueing map and $\Si_c$ the glueing of $\Si$ at the ordered pair
$((p_-,v_-),(p_+,v_+))$ with respect to $c$. Let $L_c \subset
H_1(\Si_c,\Z)$ be the Lagrangian subspace constructed above from
$L$. Then the marked surface ${\Sib}_c = (\Si_c,P,V,L_c)$ is
defined to be the {\em glueing} of ${\Sib}$ at the ordered pair
$((p_-,v_-),(p_+,v_+))$ with respect to $c$.
\end{definition}

We observe that glueing also extends to morphisms of marked
surfaces which preserves the ordered pair $((p_-,v_-),(p_+,v_+))$,
by using glueing maps which are compatible with the morphism in
question.

We can now give the axioms for a 2 dimensional modular functor.

\begin{definition} \label{DefLS}

A {\em label set\/} $\L$ is a finite set furnished with an
involution $\l \mapsto \hat \l$ and a trivial element $0$ such
that $\hat 0 = 0$.
\end{definition}

\begin{definition} \label{lmsurface}

Let $\L$ be a label set. The category of {\em $\L$-labeled marked
surfaces\/} consists of marked surfaces with an element of $\L$
assigned to each of the marked point and morphisms of labeled
marked surfaces are required to preserve the labelings. An
assignment of elements of $\L$ to the marked points of ${\Sib}$ is
called a labelling of ${\Sib}$ and we denote the labeled marked
surface by $({\Sib},\l)$, where $\l$ is the labelling.
\end{definition}

\begin{remark}{\em
The operation of disjoint union clearly extends to labeled marked
surfaces. When we extend the operation of orientation reversal to
labeled marked surfaces, we also apply the involution $\hat \cdot$
to all the labels. }\end{remark}

\begin{definition} \label{DefMF}
A {\em modular functor\/} based on the label set $\L$ is a functor
$V$ from the category of labeled marked surfaces to the category
of finite dimensional complex vector spaces satisfying the axioms
MF1 to MF5 below.
\end{definition}

\subsubsection*{MF1} {\it Disjoint union axiom\/}: The operation of disjoint
union of labeled marked surfaces is taken to the operation of
tensor product, i.e. for any pair of labeled marked surfaces there
is an isomorphism $$ V(({\Sib}_1,\l_1) \sqcup ({\Sib}_2,\l_2)) )
\cong V({\Sib}_1,\l_1) \otimes V({\Sib}_2,\l_2). $$ The
identification is associative.

\subsubsection*{MF2} {\it Glueing axiom\/}: Let ${\Sib} $ and ${\Sib}_c$ be
marked surfaces such that ${\Sib}_c$ is obtained from ${\Sib} $ by
glueing at an ordered pair of points and projective tangent
vectors with respect to a glueing map $c$. Then there is an
isomorphism $$ V({\Sib}_c,\lambda) \cong \bigoplus_{\m \in \L}
V({\Sib},\m,\hat \m,\l), $$ which is associative, compatible with
glueing of morphisms, disjoint unions and it is independent of the
choice of the glueing map in the obvious way (see remark
\ref{remarkglue2}). This isomorphism is called the glueing isomorphism and its inverse is called the factorization isomorphism.

\subsubsection*{MF3} {\it Empty surface axiom\/}: Let $\emptyset$ denote
the empty labeled marked surface. Then $$ \Dim V(\emptyset) = 1.
$$

\subsubsection*{MF4} {\it Once punctured sphere axiom\/}: Let $\Sib = (S^2,
\{p\},\{v\},\{0\})$ be a marked sphere with one marked point. Then $$
\Dim V(\Sib,\l) = \left\{ \begin{array}{ll} 1,\qquad &\l = 0\\
0,\qquad & \l \ne 0.\end{array}\right. $$

\subsubsection*{MF5} {\it Twice punctured sphere axiom\/}: Let $\Sib = (S^2,
\{p_1,p_2\},\{v_1,v_2\},\{0\})$ be a marked sphere with two marked
points. Then $$ \Dim V(\Sib,(\l,\mu)) = \left\{ \begin{array}{ll} 1,
\qquad &\l = \hat \mu\\ 0,\qquad &\l \ne \hat \mu.\end{array}\right.
$$

In addition to the above axioms one may has extra properties,
namely

\subsubsection*{MF-D} {\it Orientation reversal axiom\/}:
The operation of orientation reversal of labeled marked surfaces
is taken to the operation of taking the dual vector space, i.e for
any labeled marked surface $({\Sib},\l)$ there is a pairring $$
\langle \cdot,\cdot\rangle : V({\Sib},\l) \otimes V(-{\Sib},\hat
\l) \ra \bC{}, $$ compatible with disjoint unions, gluings and
orientation reversals (in the sense that the induced isomorphisms
$ V({\Sib},\l) \cong V(-{\Sib},\hat \l)^*$ and $V(-{\Sib},\hat \l)
\cong V({\Sib},\l)^*$ are adjoints).

\vskip.4cm

 and

\subsubsection*{MF-U} {\it Unitarity axiom\/}

Every vector
space $V({\Sib},\l)$ is furnished with a unitary structure
$$ ( \cdot,\cdot ) : V({\Sib},\l) \otimes \overline{V({\Sib},\l)}
\to {\mathbb C} $$ so that morphisms induces unitary
transformation. The unitary structure must be compatible with
disjoint union and glueing. If we have the orientation reversal
property, then compatibility with the unitary structure means that
we have a commutative diagram 
$$\begin{CD} V({\Sib},\l) @>>\cong>
V(-{\Sib},\hat \l)^*\\ @VV\cong V @V\cong VV\\
\overline{V({\Sib},\l)^*} @>\cong>> \overline{V(-{\Sib},\hat \l)},
\end{CD}$$
where the vertical identifications come from the unitary
structure and the horizontal from the duality.

\section{The skein theory construction of modular categories}

\label{cha:notes-part-one}
Let us briefly review Blanchet's \cite{Bl1} constructions of the
Hecke-category and its associated modular tensor categories to fix notation and normalisation. This
construction is really a generalization of the BHMV-construction \cite{BHMV2} of the
$U_q(sl_2(\bC))$-Witten-Reshetikhin-Turaev TQFT \cite{RT1,RT2,TW2} to the $U_q(sl_N(\bC))$-case.
We give a slightly more direct construction of this category and its associated modular functor,
which implements skein theoretically some of the abstract categorical
constructions presented in \cite{Bl1} and \cite{Tu}. This is done
in complete parallel to the $N=2$ case treated in \cite{BHMV2}.

Throughout we will fix integers, $N \geq 2$ and $K \geq 1$. Let $q$ be the following primitive $(N+K)$'th root of $1$ in $\bC$, $q=e^{2\pi i/(K+N)}$. We will also need the following roots of $q$,   $q^{1/2N} = e^{2\pi i/(2N(K+N))}$ and $q^{1/2} = e^{2\pi i/(2(K+N))}$.
We observe that the quantum integers
\[[j] = \frac{q^{j/2}-q^{-j/2}}{q^{1/2}-q^{-1/2}}\]
are invertible if $1\leq j < N+K$.

\subsection{The Hecke algebra and Jones-Wenzl idempotents}
\label{sec:hecke-algebra-proj}

Let $B_n$ be the braid group on $n$ strands.
The standard generators of $\sigma_i\in B_n$, $i=1,\dots,n-1$ are
given by the braids on $n$ strands where the $i$'th strand is
crossing over the $(i+1)$'th strand
  \begin{equation*}
  \sigma_i=\hbox{\knot{iileftover}}
  \end{equation*}

The Hecke algebra
$H_n$ is the following quotient of the group ring of $B_n$.
\begin{definition} The Hecke algebra $H_n$ is
$$
  H_n = \bC [B_n]/ \langle q^{1/2N}\sigma_i - q^{-1/2N} \sigma_i^{-1}
  =(q^{1/2}-q^{-1/2})\Id | i=1,\ldots n-1\rangle.
$$
\end{definition}

The Jones-Wenzl idempotent of $H_n$ are given explicitly as
follows:
  \begin{align*}
    g_n&=\frac{1}{[n]!}q^{-n(n-1)/4}\sum_{\pi\in S_n}
    (-q)^{(1-N)\ell(\pi)/2N}w_\pi\\
    f_n&=\frac{1}{[n]!}q^{n(n-1)/4}\sum_{\pi\in S_n} q^{-(1+N)\ell(\pi)/2N}w_\pi
  \end{align*}
where $w_\pi$ is the positive braid associated to $\pi$ and
$\ell(\pi)$ is the length of $\pi$ and the quantum factorial
$$
  [n]!= \prod_{j=1}^n [\,j\,]
$$
is
assumed to be invertible.

Following further Jones and Wenzl, we introduce the idempotents $e_i\in H_n$, $i=1,2, \ldots, n-1$ given by
$$e_i = \frac{q-q^{(N+1)/2N}\sigma_i}{q+1}.$$

These idempotents satisfies the relations (H1) and (H2) from \cite{Wzl1} and they clearly also generate $H_n$.
One can define a $*$-structure on $H_n$ by the assignment $e_i^* = e_i$, $i=1,2, \ldots, n-1$.

\subsection{Skein theory}
\label{sec:introduction}

Let $D^2$ be the unit disc in the complex plane.

\begin{definition}\label{framedsetofpoints}
    A framed set of points $\ell$ is a finite set of points $P\subset
    D^2-\partial D^2$ together with oriented directions
$$v_P \in  P(T_pD^2)^{\times P}$$ and signs $\epsilon_p$
attached to each point $p\in P$. For a framed set of points $\ell$
in $D^2$, we denote by $-\ell$ the same framed set of points, but with all
signs negated.
\end{definition}

For the definition of an oriented ribbon graph $R$ in $D^2\times [0,1]$ with
$$\partial R = \ell_0\times\{0\} \coprod -\ell_1\times\{0\}$$
 for two sets of framed points in $D^2$, see pp. 31 -- 35 in Turaev's book \cite{Tu}. We
just recall here that the signs at the boundary indicates the
direction of the band, positive for outgoing and negative for
ingoing.

We are only interested in the equivalence class of the ribbon graphs
in $D^2\times [0,1]$ up to the action of orientation preserving diffeomorphisms of
$D^2\times [0,1]$, which are the identity on the boundary and isotopic to the
identity among such.

\begin{definition}\label{arg}
A ribbon graphs $R$ in $D^2\times [0,1]$ is called special if it only contains coupons of the following
type:
 
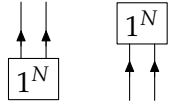
\begin{center}\label{frHomfly}

  \includegraphics{andersen-16.mps}\qquad \includegraphics{andersen-17.mps}
  \end{center}
which has $N$ incoming or $N$ outgoing bands.
\end{definition}

The label $1^N$ on these coupons is as such immaterial, but will be justified by the
relations on them introduced below.

Let  $\h(D^2\times [0,1],\ell_0,\ell_1)$ be the free complex vector space generated by special ribbon
graphs in $D^2\times [0,1]$, whose boundary is $\ell_0\times\{0\} \coprod -\ell_1\times\{0\} $, modulo the following local
relations \vskip.3cm
\begin{gather}\label{skeinrel}
q^{1/2N}\hbox{\knot{leftover}}  -q^{-1/2N} \hbox{\knot{rightover}}
=\big(q^{1/2}-q^{-1/2}\big) \hbox{\knot{nocross}}
\nonumber\\
\hbox{\knot{curvunder}}
=q^{(N^2-1)/2N}\hbox{\knot{onestraight}} 
\qquad\qquad
\hbox{\knot{curvover}}
=q^{(1-N^2)/2N}\hbox{\knot{onestraight}}
\\
L\cup \hbox{\knot{zero}} 
=[N] L \nonumber
\end{gather}
plus the two coupon relations:
\begin{center}
  \label{cr}
  \includegraphics{coupon.mps}
\end{center}

We observe
that if the algebraic number of points in $\ell_0\coprod (-\ell_1)$ is not a non-zero
multiple of $N$, then all coupons of any special ribbon graph in
$(D^2\times [0,1],\ell_0,\ell_1)$ can be reduced away by the first of the coupon relations,
so certainly $\h(D^2\times [0,1],\ell_0,\ell_1)$ is in this case generated by ribbon links.
By Proposition 1.11 in \cite{Bl1} it also follows that no extra
relations is obtained from the second coupon relation, so that
$\h(D^2\times [0,1],\ell_0,\ell_1)$ is in this case isomorphic to the Homfly skein module
of $D^2\times [0,1],\ell_0,\ell_1)$ in the $sl_N(\bC)$-specialization we are considering here.
We further observe that if $\ell_0 = \ell_1$ then framed braids spans $\h(\Sigma\times
[0,1],-\ell_\Sigma,\ell_\Sigma)$.

We have the following fundamental theorem.
\begin{theorem} The framed version of the Homfly polynomial
which by definition satisfies the relations (\ref{frHomfly}) induces
by evaluation an isomorphism
  \begin{equation*}
    \morf{\langle\cdot \rangle}{\h(\Sp^3)}{\bC}
  \end{equation*}
where $\h(\Sp^3)$ is the Homfly skein module of $\Sp^3$.
\end{theorem}

The Homfly polynomial was first introduced in \cite{Homfly}. 

\subsection{The Hecke category $H$}
\label{sec:hecke-category-h}
The Hecke category $H$ is defined as follows.
The  objects are pairs $\alpha=(D^2,\ell)$, where $\ell$ is a framed set
of points
in the interior of the $2$-disk $D^2$. The morphisms  $\Hom(\alpha,\beta)=
H(\alpha,\beta)$,
between two objects  $\alpha=(D^2,\ell_0)$ and
$\beta=(D^2,\ell_1)$ are
  \begin{equation*}
    H(\alpha,\beta)=\h(D^2\times [0,1],\ell_0\times \{0\} \amalg
    -\ell_1\times\{1\}).
  \end{equation*}
There is a trace $\Tr_\alpha^{N,K}$ on $H_\alpha=H(\alpha,\alpha)$ given by
  \begin{equation*}
    \morf{\Tr_\alpha^{N,K}\ \ }{H_\alpha}{\h(D^2\times\Sp^1)}
    \longrightarrow\h(\Sp^3)\longrightarrow \bC
  \end{equation*}
The first map is obtained by glueing the bottom and top disk of $D^2
\times [0,1]$. The second map is induced by the standard inclusion
of $D^2 \times S^1$ into $S^3$ and the last is induced by the framed
Homfly polynomial.  Define
$$
\langle\cdot, \cdot\rangle : H(\alpha,\beta) \times H(\beta,\alpha)
\ra \bC
$$
by
$$
\langle f, g\rangle = \Tr_\alpha^{N,K}(fg).
$$
The definition of tensor product, braiding, twist and duality is straight forward and explained in detail in \cite{Bl1}, where the following proposition is established.

\begin{proposition}\label{Ribcat}
  The category  $H$ with the above structure is a ribbon category.
\end{proposition}

The object consisting of the points $-(n-1)/n, -(n-3)/n, \ldots,
(n-1)/n$ framed along the positive real axis we denote simply $n$. We
observe that $H_n$ is naturally isomorphic to $H(n,n)$. The trace
$\Tr_n^{N,K}$ defined above induces therefore a trace on~$H_n$.

\subsection{Young symmetrizers}
\label{sec:young-symmetries}
Let us now for each Young diagram define an object in the Hecke
category and the Aiston-Morton realization of the corresponding Young symmetrizer following \cite{Bl1}, \cite{Aiston} and \cite{AM}.

To a \emph{partition} of $n$,
  $\lambda=(\lambda_1\geq\dots\geq\lambda_p\geq1)$,
  $n=\lambda_1+\dots+\lambda_p$ there is the Young diagram of size
  $\abs{\lambda}=n$:
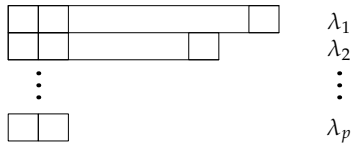
  \begin{center}
    \includegraphics[scale=0.8]{andersen-19.mps}
  \end{center}
For a Young diagram $\lambda$ we use the usual notation $\lambda^\vee$
to denote the Young diagram obtained from $\lambda$ by interchanging the rows
and columns.

Let now $\lambda$ be a Young diagram, $\abs{\lambda}=n$.
Let $\bo$ be the object in $H$ obtained by ``putting $\lambda$
  over $D^2$'', i.e. put a point at $\frac{(k+i l)}{(n+1)}$ if
  $\lambda$ has a cell at $(n,l)$, where we index by (row, column).
Let $F_\lambda\in H_{\bo}$ be obtained by putting
 $[\lambda_i]! f_{\lambda_i}$ along row $i$ in $\lambda$, $i=1,\dots,p$ and
 $G_\lambda\in H_{\bo}$ be obtained by putting
  $[\lambda_j^\vee]!g_{\lambda_j^\vee}$ along colom $j$, $j=1,\dots,p^\vee$.
Then $\tilde{y}_\lambda = F_\lambda G_\lambda$ is a
quasi-idempotent, since by Proposition 1.6 in \cite{Bl1}
$$\tilde{y}_\lambda^2 = [hl(\lambda)] \tilde{y}_\lambda,$$ where
$hl(\lambda)$ is the hook-length of
    $\lambda$. So when $[hl(\lambda)]$ is non-zero we define
     $$y_\lambda=[hl(\lambda)]^{-1}\tilde{y}_\lambda,$$ which is an
     idempotent.
By Proposition 1.8 in \cite{Bl1}, we have that
$$\mu\neq \lambda \Rightarrow y_\lambda H(\bo,\bo[\mu])y_\mu=0$$ and
$$y_\lambda H_{\bo} y_\lambda= \bC y_\lambda.$$

\subsection{$\Gamma_{N,K}$-completed Hecke category $H^{\Gamma_{N,K}}$}
\label{sec:C-compHecke-cat}

Consider the following subset of Young diagrams
\begin{equation*}
\Gamma_{N,K} = \left\{(\lambda_1,\dots,\lambda_p)\mid
    \lambda_1\leq K,p<N\right\}.
    \end{equation*}

We observe that $[hl(\lambda)]$ is non-zero for all $\lambda \in \Gamma_{N,K}$.
We now introduce the $\Gamma_{N,K}$-completed Hecke category $H^{\Gamma_{N,K}}$, whose objects are
    triples
     $\alpha=(D^2,\ell,\lambda)$ where
    $\lambda=\big(\lambda^{(1)},\dots,\lambda^{(m)}\big)
    $, $\lambda^{(i)}\in \Gamma_{N,K}$, and $\ell = (\ell_1,\ldots,\ell_m)$ being a
    framed set of points in the interior of $D^2$.
 We have an expansion operation $E$ which maps objects of $H^{\Gamma_{N,K}}$ to objects of $H$.
 For an object of
 $H^{\Gamma_{N,K}}$, we let $E(\alpha)=\big(D^2,E(\ell)\big)$ be the object in
    $H$, where $E(\alpha)$ is obtained by embed $\bo[\lambda^{(i)}]$  in a
        neighborhood of $ \ell_i$    according to the tangent vector of $ \ell_i$.
 Then $\pi_{\alpha} = y_{\lambda^{(1)}}\otimes\dots\otimes y_{\lambda^{(n)}}$ defines
    an idempotent in $H_{E(\alpha)}$ and we let
     $H^{\Gamma_{N,K}}(\alpha,\beta)=\pi_{\alpha}
     H\big(E(\alpha),E(\beta)\big)\pi_\beta$. By simply associating
the Young diagram $\Box$ to all points, we have a natural inclusion of the
objects of $H$ into the objects of $H^{\Gamma_{N,K}}$. Moreover, for such
objects, $E$ is just the identity and the hom-spaces between such
are identical for $H$ and $H^{\Gamma_{N,K}}$.
For any $\lambda\in\Gamma_{N,K}$, we simply write $\lambda$ for the object in
$H^{\Gamma_{N,K}}$ given by $\big(D^2,(0,v,+1),\lambda\big)$, where $v$ is the
direction of the positive real line through $0$. The category
$H^{\Gamma_{N,K}}$ inherits the structure of a ribbon category from $H$.

\subsection{The modular category $\hsu$}
\label{sec:modular-category-hsu}

We define the category $\hsu$ to have the objects of $H^{\Gamma_{N,K}}$ and
the morphims given by
  \begin{equation*}
    \hsu(\alpha,\beta)=\frac{H^{\Gamma_{N,K}}(\alpha,\beta)}{\N^{\Gamma_{N,K}}(\alpha,\beta)}
  \end{equation*}
 where
 $$\N^{\Gamma_{N,K}}(\alpha,\beta) = \Big\{f\in
    H^{\Gamma_{N,K}}(\alpha,\beta) \mid \inner{f,g}=0, \ \forall g\in
    H^{\Gamma_{N,K}}(\beta,\alpha)\Big\}.
   $$

This construction on the morphisms is called ``purification'', where one removes
  ``negligible'' morphisms $\N^{\Gamma_{N,K}}(\alpha,\beta)$. We observe that
  $H^{\Gamma_{N,K}}(\alpha,\beta)$ is a sub-quotient of
  $H(E(\alpha),E(\beta))$.

\noindent Let $l^m$ be the Young diagram with $l$ columns containing $m$
cells. The object 
$$1\otimes \ldots \otimes 1$$ 
($l$ factors) in $\hsu$
will simply be denoted $l$. We further use the notation 
$$\N_n^{N,K} = \N^{\Gamma_{N,K}}(n,n).$$

For a Young diagram $\lambda$
in $\overline{\Gamma}_{N,K}$ we define
  $\lambda^\dagger \in\Gamma_{N,K}$ to be the Young diagram
  obtained from the skew-diagram $(\lambda_1^N)/\lambda$ by rotation
  as indicated in the figure below. 
  
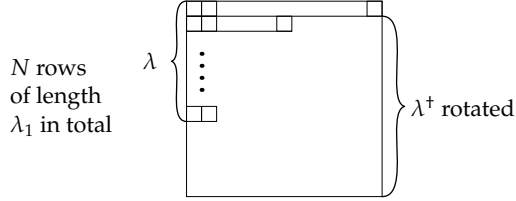
\begin{center}
    \includegraphics[scale=0.8]{andersen15aa.mps}
  \end{center}

According to Theorem 2.11 in \cite{Bl1} we have that.

  \begin{theorem}\label{BlModular}
    The category $\hsu$  with simple objects
$\lambda\in\Gamma_{N,K}$ and duality involution
  $\lambda \mapsto\lambda^\dagger$ is modular.
  \end{theorem}

Let $n\in \bN$ and consider the decomposition 
$$ \hsu(n,n) = \bigoplus_{\lambda\in \Gamma_{N,K}} \hsu(n, \lambda)\otimes \hsu(\lambda, n),$$
which we have, since $\hsu$ is a modular tensor category. Let 
$$\Gamma_{N,K}^n = \{ \lambda \in \Gamma_{N,K} \mid n\geq |\lambda| \mbox{ and } N | (n-|\lambda|)\}.$$
Let $\lambda \in \Gamma_{N,K}^n$. Then $\hsu(n,\lambda) \neq 0$ and we let $z^{(n)}_\lambda\in \hsu(n, \lambda)\otimes \hsu(\lambda, n)$ be such that
\begin{equation}\label{1}
1_n = \sum_{\lambda\in \Gamma^n_{N,K}}  z^{(n)}_\lambda.
\end{equation}
We recall that there is a one to one correspondence between irreducible representations of a finite dimensional algebra over the complex numbers and its minimal central idempotents. Each subalgebra $ \hsu(n, \lambda)\otimes \hsu(\lambda, n)$ is the full matrix algebra, since $\hsu(\lambda, n)$ is the dual of $ \hsu(n, \lambda)$ (see Corollary 4.3.2 in \cite{Tu}), so we conclude that the $z^{(n)}_\lambda$, $\lambda\in \Gamma^n_{N,K}$ are the minimal central idempotents of $\hsu(n,n)$ and $\hsu(n, \lambda)$, $\lambda\in \Gamma^n_{N,K}$, are the irreducible modules of $\hsu(n,n)$. We have thus proved the following proposition.

\begin{proposition}\label{Skeinreptheory}

We have that $\hsu(n,\lambda)$, $\lambda\in \Gamma^n_{N,K}$, are the irreducible modules of $\hsu(n,n)$. Moreover  $\hsu(n,\lambda)$ corresponds to the minimal central idempotent $z^{(n)}_\lambda$ for all $\lambda\in \Gamma^n_{N,K}$.
\end{proposition}

\begin{definition}
Suppose $n\in \bN$ and let $\lambda \in \Gamma^n_{N,K}$.
Then we define 
$$\Gamma_{N,K}^{n-1,\lambda} = \left\{ \lambda' \in \Gamma_{N,K} \mid ( \lambda' < \lambda \mbox{ and } |\lambda/\lambda'| = 1) \mbox{ or } (\lambda<\lambda' \mbox{ and }  \lambda'/\lambda = 1^{N-1}) \right\}.$$
\end{definition}

We observe that $\Gamma_{N,K}^{n-1,\lambda} \subset \Gamma_{N,K}^{n-1} $ and further that $|\Gamma_{N,K}^{n-1,\lambda}|\geq 1$ and if $n >  |\lambda| $ and $\lambda \neq \emptyset$ then  $|\Gamma_{N,K}^{n-1,\lambda}|\geq 2$. Further if $|\Gamma_{N,K}^{n-1,\lambda}| = 1$ then $\lambda = l^m$ for positive $l$ and $m$ such that $n= |\lambda|$.

\begin{proposition}

Suppose $\lambda \in \Gamma_{N,K}$ and $n\in \bN$ such that $N$ divides $n-|\lambda|$ and let $\lambda' \in \Gamma_{N,K}^{n-1,\lambda} $. Then
the $H_{n-1}$-module $\hsu(n-1,\lambda')$ is isomorphic to a submodule of the $H_{n-1}$-module $\hsu(n, \lambda)$.

\end{proposition}

\begin{proof}
  Let us consider first the case where $\lambda'<\lambda$. Let we
  observe that
$$\hsu(\lambda'\otimes 1, \lambda)\neq 0$$
by Corollary 1.10 in \cite{Bl1}, so we choose $e\in
\hsu(\lambda'\otimes 1, \lambda)- \{0\}$. By the same Corollary we see
that there then must exist $e'\in \hsu(\lambda, \lambda'\otimes 1)$,
such that $ e \circ e' = 1_\lambda$. The morphism
$$ f\mapsto e \circ  (f\otimes 1) \in \hsu(n, \lambda)$$
for $f\in \hsu(n-1,\lambda')$ must be non-zero, since it is possible
to choose $f$ such that there exists $f' \in \hsu(\lambda', n-1) $
which satisfies $f\circ f' = 1_{\lambda'}$. But then
$$e \circ (f \otimes 1)\circ (f'\otimes 1) e' = 1_\lambda$$ by Corollary 1.10 in \cite{Bl1}. Since $ \hsu(n-1,\lambda')$ is an irreducible $H_{n-1}$-module this morphism must be injective and hence an isomorphism onto its image.
Now we consider the case where $\lambda<\lambda',$ $n> |\lambda|$ and
$\lambda'/\lambda = 1^{N-1}$. First we fix an isomorphism $e\in
\hsu(1^N\otimes \lambda, \lambda)$. Now let $\tilde \lambda\in
\overline \Gamma_{N,K}$ be the young diagram obtain from $\lambda$ by
putting $1^N$ to the immediate left of $\lambda$. Again by Corollary
1.10 of \cite{Bl1}, we have that there exist
$$ g_1 \in \hsu(\lambda'\otimes 1, \tilde \lambda) - \{0\} \mbox{ and } g_2 \in \hsu( \tilde \lambda,1^N\otimes \lambda) - \{0\},$$
for which there exists
$$ g'_1 \in \hsu( \tilde \lambda, \lambda'\otimes 1) - \{0\} \mbox{ and } g'_2 \in \hsu(1^N\otimes \lambda, \tilde \lambda) - \{0\},$$
such that
$$ g_1 g_1' = 1_{\tilde\lambda}\mbox{, } g_2' g_2 = 1_{\tilde\lambda}.$$
Let $g = g_2g_1.$ Now consider the morphism
$$ f\mapsto e \circ g \circ (f\otimes 1) \in \hsu(n, \lambda),$$
for $f$ in $\hsu(n-1+N, \lambda')$. Again by Corollary 1.10 in
\cite{Bl1} we may choose $f$ such that there exist $f'\in
\hsu(\lambda', n-1+N)$, such that $f \circ f' = 1_ {\lambda'}$. But
then we see that
$$ g_2' \circ e^{-1} \circ e  \circ g  \circ (f\otimes 1) \circ (f'\otimes 1)  \circ g_1' = 1_{\tilde \lambda}.$$
hence the above morphism is non-zero and therefore again an
isomorphism onto its image.  
\end{proof}

\section{The Witten-Reshetikhin-Turaev modular functor via skein theory} \label{RTmf}

Let us briefly recall Turaev's construction of a modular functor from a modular tensor category in \cite{Tu}. Since we are interested in the modular tensor category $\hsu$, we
will at the same time apply it to Blanchet's category, the construction of which we reviewed in section \ref{cha:notes-part-one}.
Let $\Sib = (\Si, P, V, L)$ be a marked surface and
$\lambda$ a labeling of it by labels from the $\Gamma_{N,K}$, hence $\lambda : P \ra \Gamma_{N,K}$. Let  $\Sib_0 = (\Si_0, P_0, V_0, L_0)$ be the standard surface (see section 1.2 Chapter IV in \cite{Tu}) of the same type as $\Sib$ . Let the genus of $\Si$ (or equivalently of $\Sigma_0$) be $g$.
We recall the following definition

\begin{definition}
The modular functor $\V^{SU(N)}_K$ obtained by applying the Reshetikhin-Turaev construction to the modular tensor category $\hsu$ associates to any labeled marked surface $(\Sib,\lambda)$ the vector space 
$\V^{SU(N)}_K(\Sib,\lambda)$ which is uniquely determined by the following property. For any morphism $\phi : \Sib_0 \ra \Sib$ of marked surfaces, there is a unique isomorphism
$$
\V^{SU(N)}_K(\phi) : \bigoplus_{\mu \in \Gamma_{N,K}^{\times g}}\hsu\left(\lambda_1\otimes \ldots \otimes \lambda_p \otimes \left(\bigotimes_{i=1}^g \mu_i\otimes \mu_i^\dagger\right),0\right)\ra \V^{SU(N)}_K(\Sib,\lambda)
$$
where $(\lambda_1,\ldots,\lambda_p)$ is the labeling of $P_0$ induced from $\lambda$ via $\phi$, such that if $\phi' : \Sib_0 \ra \Sib$ is another parametrization, then we get the formula
$$
\V^{SU(N)}_K((\phi'))^{-1}\circ \V^{SU(N)}_K(\phi) = \varphi(\phi',\phi)
$$
where $\varphi(\phi',\phi)$ is defined in section 6.3 in Chapter IV of \cite{Tu}.
\end{definition}
We remark that $\varphi(\phi',\phi)$ is constructed by producing a ribbon graph presentation of the mapping cylinder of $(\phi')^{-1}\circ \phi$ and then computing the TQFT morphism determined by $\hsu$ for this ribbon graph . In case $g=0$, we observe that this mapping cylinder is just determined by a spherical braid on $p$ strands connecting $P_0$ to it self. For further details on this please see section 2 of Chapter IV in \cite{Tu}.

The construction of the gluing morphism for the modular functor is described in Section 4 of Chapter V in \cite{Tu}. This morphism is also constructed by producing an explicit $3$-dimensional cobordism between the unglued and the glued marked surface. This cobordism is described in detail in section 4.5 of Chapter V in \cite{Tu}.   As it is explained in section 5.10 of Chapter 10 in \cite{Tu}, the glueing isomorphism for standard genus zero surfaces is given by the following isomorphism
\begin{eqnarray} \label{facc}
\mbox{  }&H^{SU(N)}_K(\lambda_1\otimes \ldots \otimes \lambda_{p+q},0) \cong \\ 
& \bigoplus_{\mu \in \Gamma_{N,K}} H^{SU(N)}_K(\lambda_1\otimes \ldots \otimes \lambda_{p}\otimes \mu^\dagger,0) \otimes H^{SU(N)}_K(\mu \otimes\lambda_{p+1}\otimes \ldots \otimes \lambda_{p+q},0). \nonumber
\end{eqnarray}

Suppose $\alpha$ is an object in $H^{SU(N)}_K$, then by including it in $\Sigma$ we get the structure of a labeled marked surface $\Sib_\alpha$. Let $\lambda_\alpha$ be the corresponding labelling of $\Sib_\alpha$. In particular for the object $n$, we have $\Sib_n$, which we will choose as the standard labeled marked surface of genus $0$ with $n$ marked points and $\lambda_n$ will denote the box-labelling of all marked points in $n$. We remark that the subindex here refers to the object $n$.

\begin{lemma}\label{obmf}
For any object $\alpha$ of $H^{SU(N)}_K$ we have a canonical isomorphism
$$\Phi_\alpha : H^{SU(N)}_K(\alpha,0) \ra \V^{SU(N)}_K(\Sib_\alpha,\lambda_\alpha)$$
such that for all $T_{\alpha',\alpha}$ representing an element of $H^{SU(N)}_K(\alpha',\alpha)$, where $\alpha'=\lambda_1\otimes \ldots \otimes \lambda_n$ for some $n$, we get a commutative diagram
$$
\begin{CD} \hsu(\alpha,0)@>\Phi_\alpha>>  \V^{SU(N)}_K(\Sib_\alpha, \lambda_\alpha)\\ 
@VV T_{\alpha,\alpha'} V @V \V^{SU(N)}_K(\phi_{T_{\alpha,\alpha}}) VV\\
 \hsu(\alpha',0)@>\Id>>  \V^{SU(N)}_K(\Sib_{\alpha'}, \lambda_{\alpha'})
\end{CD}
$$
where $\phi_{T_{\alpha,\alpha'}}$ is the diffeomorphism of the disc, determined by the braid $T_{\alpha',\alpha}$, such that $\phi_{T_{\alpha,\alpha'}}(\alpha) = \alpha'$, extended to a morphism of labelled marked surfaces from 
$(\Sib_\alpha, \lambda_\alpha)$ to $(\Sib_{\alpha'}, \lambda_{\alpha'})$.
\end{lemma}

\section{The representation theory of the Hecke algebra}
\label{sec:rep-Hecke}

Following Jones \cite{J-4} and Wenzl \cite{Wzl1} we will now consider the tower of algebras $H_1 \subset H_2\subset \ldots$ and consider the towers of representations corresponding to Markov traces on this tower, first constructed by Wenzl in \cite{Wzl1}. We will see that the modular category $\hsu$ constructed above yields another realization of Wenzl's representations via the GNS construction. These representations are central to Wenzl's construction of subfactors from the Hecke-algebras in \cite{Wzl1}.

Let $H_\infty$ by the inductive limit of the algebras 
$$H_1\subset H_2 \subset \ldots.$$
A trace $\tr$ on $H_\infty$ is a linear functional $\tr : H_\infty \ra \bC$, such that $\tr(xy) = \tr(yx)$ for all $x,y \in H_\infty$ and such that $\tr(1) = 1$.

\begin{definition}

A trace $\tr$ on $H_\infty$ is a Markov trace if there is an $\eta \in \bC$ such that 
$$\tr(xe_n) = \eta \tr(x)$$
for all $x\in H_n$ and all $n\in \bN$.

\end{definition}

\begin{lemma}
A Markov trace is uniquely determined by $\eta=\tr(e_1)$.
\end{lemma}

For an argument for this lemma see \cite{J-4}. In \cite{J-4} Jones also proved the theorem, that there exist a Markov trace for any $\eta\in \bC$. Let now
$$\eta_{N,K} = (q-q^N)/(1+q)(1-q^N),$$
and let $\tr^{N,K}$ be the corresponding Markov trace.

\begin{theorem}[Wenzl]\label{Wenzl}
The GNS construction applied to the state $\tr^{N,K}$ on the *-algebra $H_\infty$ gives a representation $\pi^{N,K+N}$ of $H_\infty$ with the property\footnote{Please see the discussion following this theorem for the proper interpretation of this property.} that 
\begin{equation}\pi^{N,K+N}\mid_{H_n} \cong \bigoplus_{\lambda \in \Lambda_n^{N,K+N} }\pi_\lambda^{N,K+N},\label{GNSW}
\end{equation}
where $ \Lambda_n^{N,K+N} $ is the set of Young diagrams defined in Definition 2.4 in \cite{Wzl1} and the $\pi_\lambda^{N,K+N}$'s are Wenzl's representations
$$
\pi_\lambda^{N,K+N} : H_n \ra B(V^{N,K+N}_\lambda)
$$
constructed right after Definition 2.4 in \cite{Wzl1}.
\end{theorem}

This theorem is one half of Theorem 3.6(b) in \cite{Wzl1}. The other half states that it is precisely these Markov traces which factors through a $C^*$-representation. 

Let us now briefly recall the GNS construction: Define 
$$\N_\infty^{N,K} = \{x \in H_\infty \mid \tr^{N,K}(x^*x) = 0\},$$
and let $\tilde H^{N,K}_\infty = H_\infty/\N_\infty^{N,K}$.
Then $\pi^{N,K+N} : H_\infty \ra B(\tilde H_\infty^{N,K})$ is just given by 
the left action of $H_\infty$ on $\tilde H_\infty^{N,K}$. Here $\tilde H_\infty^{N,K}$ gets a pre-Hilbert space structure induced from the trace by the formula
$$(x,y)_{N,K} = \tr^{N,K}(y^*x).$$
The Cauchy-Schwarz inequality shows that this is indeed a well-defined inner product on $\tilde H_\infty^{N,K}$.
The identity $1\in H_\infty$ projects in $\tilde H_\infty^{N,K}$ to the required cyclic vector for 
$$\pi^{N,K} : H_\infty \ra B(\tilde H_\infty^{N,K}).$$
The isomorphism in Theorem \ref{Wenzl} should be understood in the following sense. Let $\tilde H^{N,K}_n$ be the image of $H_n$ in $\tilde H_\infty^{N,K}$. Then the action of $H_n$ on $\tilde H_\infty^{N,K}$ preserves $\tilde H_n^{N,K}$ and it is this representation of $H_n$ on $\tilde H_n^{N,K}$ which is isomorphism to the direct sum of Wenzl's representations given in the right hand side of (\ref{GNSW}). We will use the notation
$$ \pi_n^{N,K+N} : H_n \ra B(\tilde H_n^{N,K})$$
for this representation.

We now introduce the following trace on $H_n$
$$\Tr^{N,K}(x) = \Tr^{N,K}_n(x)/\Tr^{N,K}_n(\id),$$
for all $x\in H_n.$ We observe that $\Tr^{N,K}(x)  = \Tr^{N,K}(x\otimes 1) $ for all $x\in H_n$, hence $\Tr^{N,K}$ is a well-defined trace on $H_\infty$.

\begin{theorem}
The traces $\Tr^{N,K}$ and $\tr^{N,K}$ coincide.

\end{theorem}

\begin{proof}
  From the skein theory it is clear that $\Tr^{N,K}$ is a Markov trace
  and an explicitly simple computation show that $\Tr^{N,K}(e_1) =
  \eta_{N,K}$.
\end{proof}

\begin{theorem}
For any $n$, we have that 
$$\N_n^{N,K} = \{x\in H_n \mid \Tr^{N,K}(x^*x) =0\}.$$
\end{theorem}

\begin{proof}
  It is clear that the left hand side is contained in the right hand
  side.  The other inclusion follows from the observation that
  $\Tr^{N,K}(yx) = \Tr^{N,K}((y^*)^*x) = 0$ for all $x \in
  \N_\infty^{N,K}$ by Cauchy-Schwarz.
\end{proof}

By using similar easy arguments one finds that $\N^{N,K}_\infty$ is two-sided ideal in $H_\infty$, from which we conclude that $\pi^{N,K+N}$ factors to a well-defined representation $\pi^{N,K+N} : H_\infty/\N^{N;K}_\infty \ra B(\tilde H_\infty^{N,K}).$ A simple argument shows that this representation is injective.
Moreover, if we consider $\pi_n^{N,K+N}$, it clearly also factors through $H_n\cap \N^{N,K}_\infty$, and the same argument again shows that also this representation is injective on $H_n/H_n\cap \N^{N,K}_\infty$ which by the above lemma is naturally isomorphic to $\hsu(n,n)$.
From this we have the following theorem as an immediate consequence.

\begin{theorem}
For all $n$ we have a canonical $H_n$ isomorphism 
$$\Psi_n^{N,K} :\hsu(n,n) \stackrel{\cong}{\ra} \pi^{N,K+N}_n(H_n),$$
which is compatible with inclusions
$$\begin{CD} \hsu(n,n)@>\Psi_n^{N,K} >> \pi^{N,K+N}_n(H_n)\\ 
@VVV @VVV\\
\hsu(n+1,n+1)@>\Psi_{n+1}^{N,K} >> \pi^{N,K+N}_{n+1}(H_{n+1})
\end{CD}
$$ 
\end{theorem}

We now consider the following inductive limit.
\begin{definition}
We let 
$$
\hsu(\infty,\infty) = \varinjlim_n \hsu(n,n).
$$
\end{definition}
From the above discussion we see the following Theorem.

\begin{theorem}\label{SW}
The isomorphisms $\Psi_n^{N,K}$ induces a natural $H_\infty$-isomorphism 
$$ \Psi_\infty^{N,K} : \hsu(\infty,\infty)\ra \tilde H_\infty^{N,K}.$$
\end{theorem}

Let $\lambda\in  \Gamma_{N,K}$. Let $n$ be an integer such that $n-|\lambda|$ is a non-negative multiple, say $l$, of $N$, thus $\lambda\in  \Gamma^n_{N,K}$. Let $ (\lambda)_n$ be obtained from $\lambda$ by attaching 
$l$ column of $N$ boxes to the left of $\lambda$. This provides a bijection between $\Gamma^n_{N,K}$ and Wenzl's $\Lambda_n^{N,K+N}$.

\begin{theorem}\label{DECOMPCOMP}
Let $n\in \bN$.
The canonical $H_n$ isomorphism $\Psi_n^{N,K} $ from the previous theorem takes the decomposition
$$\hsu(n,n) \cong \bigoplus_{\lambda\in \Gamma^n_{N,K}} \hsu(n,\lambda) \otimes \hsu(\lambda,n)$$
to the decomposition
$$\pi_n^{N,K+N}(H_n) = \bigoplus_{\lambda \in  \Gamma^n_{N,K}} \pi_{(\lambda)_n}^{N,K+N}(H_n),$$
for all $n$. 
Moreover for any $\lambda\in  \Gamma^n_{N,K}$ we have that $\hsu(n,\lambda)$ and $V^{N, K+N}_{(\lambda)_n}$
are isomorphic representations of $H_n$.
\end{theorem}

\begin{proof}
  Following again Wenzl, we consider the minimal central idempotents
  $\tilde z_{(\lambda)_n}$ corresponding to $
  \pi_{(\lambda)_n}^{N,K+N}(H_n) \subset \pi^{N,K+N}(H_n)$. Likewise
  we consider the minimal central idempotents $z^{(n)}_\lambda \in
  \hsu(n,n)$ corresponding to
$$\hsu(n,\lambda) \otimes \hsu(\lambda,n) \subset \hsu(n,n).$$ 
We claim that
\begin{equation}\Psi_n^{N,K}(z^{(n)}_\lambda) = \tilde
  z_{(\lambda)_n}\label{Idci}.
\end{equation}

It is clear that this formula for a given $n$ is implied by the second
half of the theorem, since there is a one to one correspondence
between minimal central idempotents and irreducible representations of
$\hsu(n,n) \cong \pi^{N,K+N}(H_n)$.

For $n=1$ the theorem is trivial true.

For $n=2$ we have the following two possible diagrams: $\lambda_{(2)}$
consisting of a row of two boxes and $\lambda_{(1,1)}$ consisting of a
column consisting of two boxes. A simple calculation with the two
minimal idempotents corresponding to the two labels $\lambda_{(2)}$
and $\lambda_{(1,1)}$ immediately gives the conclusion of the theorem
for $n=2$.

Assume that we have proved the theorem for all integers less than
$n$. We will now deduce the second half of the theorem for the integer
$n$.  Suppose $\lambda\in \Gamma^n_{N,K}$.  By Corollary 2.5 in
\cite{Wzl1} and the discussion following the proof of the Corollary,
we have that $V^{N, K+N}_{(\lambda)_n}$ is the irreducible
representation of $\pi^{N,K}(H_n) $ which corresponds to $\tilde
z_{(\lambda)_n}$. By Proposition \ref{Skeinreptheory} we have that
$\hsu(n,\lambda)$ is the irreducible representation of $\hsu(n,n)$
corresponding to $z^{(n)}_\lambda$. We now consider these to
representations as representations of $H_{n-1}$. By the explicit
construction of the representation $V^{N, K+N}_{(\lambda)_n}$ in
\cite{Wzl1} we see that it decomposes as follows under the action of
$H_{n-1}$
$$ V^{N, K+N}_{(\lambda)_n} = \bigoplus_{ \lambda' \in \Gamma^{n-1,\lambda}_{N,K} }  V^{N, K+N}_{(\lambda')_{n-1}}.$$
Then by induction we get an $H_ {n-1}$ module isomorphism
$$V^{N, K+N}_{(\lambda)_n} \cong \bigoplus_{\lambda' \in \Gamma_{N,K}^{n-1,\lambda} } \hsu(n-1,\lambda').$$
But since $V^{N, K+N}_{(\lambda)_n}$ is an irreducible
$\hsu(n,n)$-module, we know there is a unique $\mu(\lambda)\in
\Gamma^n_{N,K}$ such that
$$V^{N, K+N}_{(\lambda)_n} \cong \hsu(n,\mu(\lambda)).$$
We observe that the map $\lambda \mapsto \mu(\lambda)$ is a bijection
from $ \Gamma_{N,K}^{n}$ to it self. We also know that
$$ \hsu(n,\mu(\lambda)) > \bigoplus_{\mu' \in \Gamma_{N,K}^{n-1,\mu(\lambda) } }\hsu(n-1,\mu').$$
Hence we conclude that $\Gamma_{N,K}^{n-1,\mu(\lambda)} \subset
\Gamma_{N,K}^{n-1,\lambda}$.

\begin{claim}
  We have that $\mu(\lambda) = \lambda$ for all $\lambda \in
  \Gamma_{N,K}^{n}$.
\end{claim}

First we observe that $\mu(\emptyset) = \emptyset$ (in the case $N
\mid n$).  Second we will consider the special case where $\lambda =
l^m$ for some $l,m\in \bN$. Then $|\Gamma_{N,K}^{n-1,\lambda }|=1$ if
$n= |\lambda| $ and $|\Gamma_{N,K}^{n-1,\lambda }|= 2$ if $n>
|\lambda| $. But we must have that $|\Gamma_{N,K}^{n-1,\mu(\lambda)}|
\geq 1$ and $|\Gamma_{N,K}^{n-1,\mu(\lambda)}|= 2$ if $n>
|\mu(\lambda)| $. In the case where $|\Gamma_{N,K}^{n-1,\lambda }|=2$
and $|\Gamma_{N,K}^{n-1,\mu(\lambda)}| =1$, we have that $\mu(\lambda)
= (l')^{m'}$ for positive $l'$ and $m'$ and $n=|\mu(\lambda)|$, which
then leads to a contradiction by examining the two possibilities for
the two non-empty strict subset $\Gamma_{N,K}^{n-1,\mu(\lambda)}
\subset \Gamma_{N,K}^{n-1,\lambda}$. We must thus have that
$\Gamma_{N,K}^{n-1,\mu(\lambda)} = \Gamma_{N,K}^{n-1,\lambda}$. But
then we see immediately that in this case $\mu(\lambda) = \lambda$.

We split the rest of the proof into four cases. In each case we assume
that $\lambda \neq \mu(\lambda)$ and then we derive a
contradiction. We further recall that we can assume that $n>2$. We
will just write $\mu= \mu(\lambda)$.

\noindent {\bf Case 1.}  $n = |\lambda| = |\mu|.$ Since
$\Gamma_{N,K}^{n-1,\mu}\neq \emptyset$, we can choose $\mu' \in
\Gamma_{N,K}^{n-1,\mu}$. Then $\mu'$ is obtained from $\lambda$ by
removing one box say $b_\lambda$. Since $\lambda\neq \mu$, we see
there is unique Yong diagram $\nu$ such that $|\nu| = n+2$ and such
that $\lambda < \nu$ and $\mu < \nu$ and a box $b_\mu$ in $\nu$ such
that
\begin{equation}\label{nml}
  \nu = \lambda \cup b_\mu\mbox{, } \nu = \mu \cup b_\lambda.
\end{equation}
If $\mu$ has another box $b$ different from $b_\mu$ such that $\mu-b$
is a Yong diagram, then $\mu-b \in \Gamma^{n-1,\mu}_{N,K}$, but $\mu
-b \not \in \Gamma^{n,\lambda}_{N,K}$, hence this is a
contradiction. But then $|\Gamma^{n,\mu}_{N,K}| = 1$. This mean that
there exists $l,m \in \bN$ such that $\mu = l^m$. But then we have a
contradiction, since we are assuming that $\mu\neq \lambda$, yet we
have established in this special case for $\mu$, that $\mu = \lambda$,
since $\lambda \mapsto \mu(\lambda)$ is a bijection.

\noindent {\bf Case 2.}  $n > |\lambda|$ and $n = |\mu|.$ By counting
boxes, we see that we cannot have the situation of equation
(\ref{nml}), hence we can assume that
$$\mu - b_\mu = 1^{N-1} \cup \lambda.$$
If $\mu$ has another box $b\neq b_\mu$ such that $\mu-b$ is a Young
diagram, then there must exists a box $b'$ in $\lambda$ such that
$$\mu-b = \lambda - b'$$
which contradicts that $|\lambda| \neq |\mu|.$

\noindent {\bf Case 3.}  $n = |\lambda|$ and $n > |\mu|.$ Then we know
that $\lambda$ has a removable box $b_\lambda$ such that
$$\lambda - b_\lambda = 1^{N-1} \cup \mu,$$
which implies that $|\lambda| = N + |\mu|$. But if $\mu \neq
\emptyset$, then $\mu$ has a box $b_\mu$ which can be removed, so
$\lambda$ must have another box $b'_\lambda$ such that
$$
\lambda - b'_\lambda = \mu - b_\mu,
$$
which contradicts the above count of sizes. Hence we must have that
$\mu=\emptyset$. From this we have an immediate contradiction by the
above.

\noindent {\bf Case 4.}  $n > |\lambda| $ and $ n> |\mu|.$ We canot
have that $1^{N-1} \cup \lambda = 1^{N-1} \cup \mu$, since $\mu\neq
\lambda$. Hence $\lambda$ must have a removable box $b_\lambda$ such
that
$$
\lambda -b_\lambda = 1^{N-1} \cup \mu.
$$
Hence $|\lambda| = N + |\mu|.$ But if $\mu \neq \emptyset$ then it has
a removable box $b_\mu$ such that
$$
1^{N-1} \cup \lambda = \mu -b_\mu \implies |\lambda| = |\mu| -N
$$
or for some other removable box $b'_\lambda$ from $\lambda$
$$
\lambda - b'_\lambda = \mu -b_\mu \implies |\lambda| = |\mu|,
$$
which contradicts $|\lambda| = N + |\mu|$. So we must have that $\mu =
\emptyset$. Again this gives an immediate contradiction.  
\end{proof}

\begin{corollary}\label{cor}
There exists $H_n$-module isomorphisms
$$ \hsu(n,\lambda)  \cong V^{N,K+N}_{(\lambda)_{n}}$$
for all $\lambda \in \Gamma_{N,K}^n$.
\end{corollary}

\begin{definition}\label{swfix}
For each $\lambda \in \Gamma_{N,K}$ fix an $H_n$-module isomorphism (unique up to scale)
$$ \Psi^{N,K}_{|\lambda|,\lambda} : \hsu(|\lambda|,\lambda) \ra V^{N,K+N}_{(\lambda)_{|\lambda|}}.$$
\end{definition}

Let us now recall the following about Wenzl's construction of the representations $V^{N,K+N}_\lambda$ for each $\lambda \in \Lambda_n^{N,K+N} $. Let $T_n(\lambda)$ denote the set of generalised Young tableau's on $\lambda$ of length $n$ as defined in Definition 2.4 in \cite{Wzl1}. Then $V^{N,K+N}_\lambda$ is the complex vector space generated by $T_n(\lambda)$. Let us denote the basis vector corresponding to $t\in T_n(\lambda)$ by $w_t$. Then $w_t, t\in T_n(\lambda)$, is an orthonomal basis for $V^{N,K+N}_\lambda$. Formula (2.9) in \cite{Wzl1} gives the explicit action on this basis by the generators $e_i$ of $H_n$. Suppose that $\lambda \in  \Lambda_n^{N,K+N} $ is obtained from  $\mu \in \Lambda_{n-1}^{N,K+N} $ by adding one box, then we get an inclusion
$$\iota_{\mu,\lambda} : V^{N,K+N}_\mu \ra V^{N,K+N}_\lambda$$
of $H_{n-1}$ modules which maps $w_{t'}$ to $w_t$ for all $t\in T_n(\lambda)$ such that the diagram of $t'$ is $\mu$.
By the very definition we have an identification
$$ \pi_{(\lambda)_n}^{N,K+N}(H_n) = \End(V^{N,K+N}_\lambda).$$
Under this identification we have that Wenzl's path idempotents $p_t \in \pi_{(\lambda)_n}^{N,K+N}(H_n)$ is identified with the rank one projection onto $w_t$ for all $t\in T_n(\lambda)$. 
From this we see that the inclusion
$$  \iota_n: \pi_{n-1}^{N,K+N}(H_{n-1}) \ra \pi_n^{N,K+N}(H_{n}) $$
is given by
$$ \iota_n =  \\ \\ \bigoplus_{\mu \in  \Gamma^{n-1}_{N,K}}\\ \\ \\ \\ \bigoplus_{\{\lambda \in  \Gamma^{n}_{N,K}\mid \mu\in \Gamma_{N,K}^{n-1, \lambda}\}} \iota_{\mu,\lambda} \otimes \iota^*_{\mu,\lambda}$$
when identifying 
$$
\pi_n^{N,K+N}(H_n) = \bigoplus_{\lambda \in  \Gamma^n_{N,K}} \End(V^{N,K+N}_\lambda),
$$
for all n.

We now need to determine an $H_n$-isomorphism
$$ \Psi^{N,K}_{n,\lambda} : \hsu(n,\lambda) \ra V^{N,K+N}_{(\lambda)_n}$$
for all $\lambda \in \Gamma_{N,K}$ and integers $n = l N + |\lambda|$, where $l$ is a positive integer. These are of course unique up to scale.

We observe that we have the inclusion
$$ \hsu(N,0)\otimes \ldots\otimes \hsu(N,0)\otimes \hsu(|\lambda|, \lambda) \subset \hsu(lN+|\lambda|,\lambda)$$
and we have the inclusion
$$ V^{N,K+N}_{(\lambda)_{|\lambda|}} \subset V^{N,K+N}_{(\lambda)_{lN+|\lambda|}}$$
obtained by pre-composing the tableaux's from $T_{|\lambda|}(\lambda)$ with the tableaux of the Young diagram $N^l$ which labels the points down the first column, then down the second and so on up to the $l$'th.  The space $\hsu(N,0)$ contains a canonical non-zero element consisting of one coupon label $1^N$ attached to the object $N$ and embed in the $[-1,1]\times I \subset D^2\times I$. By tensoring this element with itself $l$ times, we get a preferred inclusion of $\hsu(|\lambda|,\lambda)$ in $\hsu(lN+|\lambda|,\lambda)$. These two inclusions allows us to fix
$$ \Psi^{N,K}_{lN+|\lambda|,\lambda} : \hsu(lN+|\lambda|,\lambda) \ra V^{N,K+N}_{(\lambda)_{lN+|\lambda|}},$$
so as to agree with the above choice for $\Psi^{N,K}_{|\lambda|,\lambda} $ under these inclusions.

Let $\lambda \in \Gamma_{N,K}$ and $t\in T_n(\lambda)$. Define $\tilde w_t \in \hsu(n,\lambda)$ by
$$ w_t  = \Psi^{N,K}_{n,\lambda} (\tilde w_t).$$

\begin{proposition}
Let $\mu,\lambda \in \Gamma_{N,K}$ such that
$$\hsu(\mu\otimes 1, \lambda) \neq 0.$$
Then there are unique 
$$\tilde w_{\mu,\lambda}\in \hsu(\mu\otimes 1, \lambda) - \{0\},$$
such that the map
$$ \hsu(n+1,\lambda) \cong \bigoplus_{\mu\in \Gamma_{N,K}}\hsu(n,\mu) \otimes \hsu(\mu\otimes 1, \lambda)$$
maps
$$ \tilde w_{\tilde t} \mapsto \tilde w_t \otimes \tilde w_{\mu,\lambda}$$
for all $t\in T_n(\mu)$ and $\tilde t \in T_{n+1}(\lambda)$ such that 
$ \tilde t' = t.$
\end{proposition}

\begin{proof}
  The direct sum decomposition
$$ \hsu(n+1,\lambda) \cong \bigoplus_{\mu\in \Gamma_{N,K}}\hsu(n,\mu) \otimes \hsu(\mu\otimes 1, \lambda)$$
of $H_n$-modules, which we get from factoring $\hsu(\nu, \lambda
\otimes 1^*)$ over the simple objects $\Gamma_{N,K}$, must be taken by
$\Psi^{N,K}_{n+1,\lambda}$ to the the direct sum decomposition
$$ V^{N,K+N}_{(\lambda)_{n+1}} \cong \bigoplus_{\mu\in \Gamma^{n,\lambda}_{N,K}} V^{N,K+N}_{(\mu)_{n}}.$$
From this we see that $\dim(\hsu(\mu\otimes 1, \lambda)) = 1$ for
exactly the $\mu \in \Gamma^{n,\lambda}_{N,K}$ and else zero and
further more there are unique elements
$$\tilde w_{\mu,\lambda}\in \hsu(\mu\otimes 1, \lambda) - \{0\}$$
with the property stated in the proposition.  
\end{proof}

The decompositions
$$\hsu(\mu\otimes (n+1), \lambda) \cong \bigoplus_{\nu\in \Gamma_{N,K}}\hsu(\mu\otimes n,\nu) \otimes \hsu(\nu\otimes 1, \lambda),$$
which we get from factoring $\hsu(\nu\otimes n, \lambda \otimes 1^*)$ over the simple objects $\Gamma_{N,K}$, together with the vectors $\tilde w_{\mu,\lambda}$ defines inductively bases for all the vector spaces $\hsu(\mu \otimes n, \lambda)$.
We observe that these bases jointly with the basis vectors $\tilde w_t $ are natural with respect to the isomorphisms
\begin{equation}\label{facdecom}
 \hsu( n+n', \lambda) \cong \bigoplus_{\mu\in \Gamma_{N,K}}\hsu( n,\mu) \otimes \hsu(\mu\otimes n', \lambda),
 \end{equation}
for all $\mu,\lambda \in \Gamma_{N,K}$ and all non-negative integers $n$ and $n'$.

\section{The genus zero isomorphism}
\label{cha:Identification}

\subsection{The label sets and the action of the Hecke algebra}
\label{sec:no-sectional-title}

Recall the construction of  the vacua modular functor $V^{\fg}_{K}$ given in \cite{AU2} for any simple Lie algebra $\mathfrak{g}$ over the complex numbers $\mathbb{C}$ and level $K$. Our main theorem states

\begin{theorem}\label{MainC}
The functor $V^{\fg}_{K}$ from the category of labeled marked
surfaces to the category finite dimensional vector spaces is a
modular functor.
\end{theorem}

We recall further that the label set for this modular functor is
\begin{equation}\label{labelset}
 P_K = \{\, \lambda \in p_+ \, | \,
   0\le(\theta, \lambda) \le K \, \}
\end{equation}
where $p_+$ is the set of dominant integral weights. Here $(\phantom{a}, \phantom{a})$ is the normalized Cartan-Killing form defined to be a constant multiple of 
the Cartan-Killing form such that
$$  
 (\theta ,\, \theta) = 2. 
$$
for the longest root $\theta$. Further the  involution 
\begin{eqnarray*}
\dagger &: &P_K \rightarrow P_K  \\
& &\lambda \mapsto \lambda^\dagger 
\end{eqnarray*}
is defined by 
\begin{equation}\label{involution}
\lambda^\dagger = - w(\lambda)
\end{equation}
where $w$ is the longest element of the Weyl group of 
the simple Lie algebra $\mathfrak{g}$.  

The modular functor constructed there for the Lie algebra $\gg= \mathfrak{sl}(N, \bC)$ at level $K$, will here be denote by $\Vdag_{N,K}$. Let us explicate the normalisation of the Cartan-Killing form and that of the longest root for this particular Lie algebra. 

Let $E_{ij}$, $1\le i, j \le N$ be the $N \times N$-matrix whose $(i,j)$-entry is
1 and all others are 0. We also let $\epsilon_i$, $i=1$, $2, \ldots, N$
be an element of the dual vector space of $\oplus_{i=1}^N \bC E_{ii}$
defined by
$$
\epsilon_i(E_{jj}) = \delta_{ij}, \quad 1\le i,j \le N.
$$
Put $r=N-1$. Also put
$$
 H_i=E_{ii}-E_{i+1, i+1}, \quad i=1,2, \ldots, r, \quad \gh = \oplus_{i=1}^r \bC H_i.
$$
Then, $\gh$ is a Cartan subalgebra of $\gg$. The dual space $\gh^*$ is written as
$$
\gh^*=\{\, \sum_{i=1}^N n_i \epsilon_i \, |\, \sum_{i=1}^N n_i = 0 \, \}.
$$
Denote by $\Delta$ the root system of $(\gg, \gh)$. A root $\alpha$
of $(\gg, \gh)$ has the form $\alpha_{ij} = \epsilon_i - \epsilon_j$, $i \ne j$,
$1 \le i, j \le N$. Choose $\alpha_i = \epsilon_i - \epsilon_{i+1}$,
$i =1,2 , \ldots, r$ as simple positive roots.

For a root $\alpha_{ij}= \alpha_i -\alpha_j$, $i\neq j$, the root space
$$
\gg_{\alpha_{ij}}= \{ \, X \in \gg\, |\, \ad (H) X = \alpha_{ij}(H)X,
\quad \forall H \in \gh\,\}
$$
is spanned by the matrix $E_{ij}$. Hence, the root space decomposition
$\gg$ is given by
$$
\gg = \gh \oplus \sum^N_{i \ne j \,  i,j =1} \bC E_{ij} .
$$
Let $(\phantom{a}, \phantom{a})$ be a constant multiple of the
Cartan-Killing form of $\gg$ defined by
$$
(X,Y) = \tr(XY), \quad X, Y \in \gg.
$$
Since the bilinear form $(\phantom{a}, \phantom{a})$ is positive definite, in the following,
we identify $\gh^*$ with $\gh$ via this bilinear form.
For any root $\alpha$ define $H_\alpha \in \gh$ in such a way that
$$
\langle \alpha, H \rangle =(H, H_\alpha)
$$
holds for any $H \in \gh$. Then for the simple roots $\alpha_i$ we have that
$H_{\alpha_i} = H_i$. The bilinear form  $(\phantom{a}, \phantom{a})$ on
the dual vector space $\gh^*$ is defined by
$$
(\alpha, \beta) = (H_\alpha, H_\beta), \quad \alpha, \beta \in \Delta.
$$
Then, for the simple roots $\alpha_i$ we have that
\begin{equation}\label{simpleroot}
(\alpha_i, \alpha_i) = 2, \quad (\alpha_i, \alpha_{i \pm 1}) = -1, \quad
(\alpha_i, \alpha_j) = 0 \quad \hbox{\rm if } |i-j|\ge 2.
\end{equation}

The longest root $\theta$ of $\gg$ is given by
\begin{equation}\label{longest}
\theta= \alpha_1+\alpha_2+\cdots + \alpha_r= \epsilon_1-\epsilon_N
\end{equation}
and we have that $(\theta, \theta) = 2$. Hence our inner product $(\phantom{a}, \phantom{a})$
is the normalized Cartan-Killing form.

The fundamental weights $\Lambda_i \in \gh^*$, $i=1,2,\ldots,r$ are defined by
\begin{equation}\label{fundamental}
\langle \Lambda_i, H_j\rangle = \delta_{ij}, \quad 1\le i,j\le r
\end{equation}
and they are given by
$$
\Lambda_i= \epsilon_1+\cdots+\epsilon_i - \frac{i}{N}\sum_{i=1}^N\epsilon_i,
\quad i=1,2, \ldots, r.
$$

\begin{lemma}
The set of dominant weights corresponding to the irreducible finite dimensional
representations of $SU(N)$ index by the Young diagrams in $\Gamma_{N,K}$ is precisely $P_K$.
\end{lemma}

\begin{proof}
 For $\lambda=(\lambda_1, \cdots, \lambda_p) \in \Gamma_{N,K}$
  we let $\Lambda$ be the corresponding dominant integral weight of
  $\mathfrak{sl}(N,\bC)$.  Then we have that
$$
\Lambda = \sum_{i=1}^{N-1} (\lambda_i - \lambda_{i+1})\Lambda_i
$$
where $\Lambda_i$ is the fundamental weight defined by
\eqref{fundamental}.  The longest root is given by \eqref{longest} and
by \eqref{simpleroot} we have that
$$
(\Lambda, \theta) = \lambda_1-\lambda_N.
$$
This gives the desired result. 
\end{proof}

We recall that the defining representation of $SU(N)$
corresponds to the Young diagram $\Box \in \Gamma_{N,K}$.

Let $\Sib_n^{\infty}$ be obtained from $\Sib_n$ by further marking $\infty$ and providing this point with the direction of the negative real axis.  Let $\mu_n = (\Box,\ldots,\Box)$ be the box labelling of the object $n$. 

We have a group homomorphism from the braid group on $n$ stands
$B_n$ to the mapping class group of $\Sib^\infty_n$, which we denote
$\Gamma_{\Sib^\infty_n}$:
$$
f : B_n \ra \Gamma_{\Sib^\infty_n}.
$$
This group homomorphism induces an algebra morphism
$$
\Upsilon_n : \bC[B_n] \ra \End(\Vdag_{N,K}(\Sib^\infty_n,\mu_n,\lambda))
$$
by the assignment
$$
\Upsilon_n(b) = \Vdag_{N,K}(f(b))
$$
for all $b\in B_n$.

For this action Kanie proved in \cite{Kanie1} the following formula. See also \cite{Ue3}.

\begin{theorem}[Kanie]\label{Skeinrelvac}
    We have the following skein relation for $\Upsilon_n$
    $$
    q^{\frac1{2N}}\Upsilon_n(\sigma_i) -
    q^{-\frac1{2N}}\Upsilon_n(\sigma_i^{-1})
    = (q^{\frac12} - q^{-\frac12}) \Id.$$

\end{theorem}

This means that $\Upsilon_n$ factors to the Hecke algebra $H_n$, and
so we see that $\Vdag_{N,K}(\Sib^\infty_n,\mu_n,\lambda)$ becomes a module over
$H_n$. Kanie constructed in \cite{Kanie1}  for all $\lambda\in \Gamma_{N,K}$ and all $n$ an explicit isomorphism 
$$\Phi^{N,N+K}_{n,\lambda}: \Vdag_{N,K}(\Sib^\infty_n,\mu_n,\lambda) \ra V^{N,K+N}_{(\lambda)_n}$$ 
by constructing a basis $v_t$, $t\in T_n(\lambda)$ of $ \Vdag_{N,K}(\Sib^\infty_n,\mu_n,\lambda)$ and defining $\Phi^{N,N+K}_{n,\lambda}$ by
$$\Phi^{N,N+K}_{n,\lambda}(v_t) = w_t$$
for all $t\in T_n(\lambda)$. He further proved the following theorem

\begin{theorem}[Kanie]\label{KanieTh}
The isomorphism $\Phi^{N,N+K}_{n,\lambda}$ is an isomorphism of $H_n$ representation from
$$
\Upsilon : H_n \ra B(\Vdag_{N,K}(\Sib^\infty_n,\mu_n,\lambda^\dagger))
$$
to 
$$
\pi_{(\lambda)_n}^{N,K+N} : H_n \ra B(V^{N,K+N}_{(\lambda)_n})
$$
for all $\lambda \in \Gamma^n_{N,K}$.
\end{theorem}

This theorem is proved in \cite{Kanie1}. See also Theorem 6.29 in \cite{Ue3}.

Kanie's basis $v_t$, $t\in T_n(\lambda)$, of $\Vdag_{N,K}(\Sib^\infty_n,\mu_n,\lambda)$ is constructed inductively such that for all $\mu,\lambda\in \Gamma_{N,K}$ and integer $n$ with the property that there is a non-negative integer $l$ solving the equation $n=lN+|\lambda|$ and $\mu\in \Gamma_{N,K}^{n,\lambda}$, he constructs a non-zero vector
$$v_{\mu,\lambda} \in \Vdag(\Sib^\infty_2, \mu,\Box, \lambda^\dagger)$$
such that if $t\in T_n(\lambda)$ such that $t'$ is a tableau for $\mu$, then
$$ v_t = v_{t'} \otimes v_{\mu,\lambda},$$
under the glueing isomorphism
$$ \Vdag(\Sib^\infty_{n}, \mu_n, \lambda^\dagger)\cong \bigoplus_{\mu \in \Gamma_{N,K}^{n,\lambda}} \Vdag(\Sib^\infty_{n-1}, \mu_{n-1}, \mu^\dagger)\otimes \Vdag(\Sib^\infty_2, \mu,\Box, \lambda^\dagger).$$

\begin{definition}\label{hsuglabel}
The isomorphisms
$$
I_{N,K}(n, \mu) : \hsu(n,\mu) \ra \Vdag_{N,K}(\Sib^\infty_n,\mu_n,\mu^\dagger) 
$$
are given by the formula
$$ I_{N,K}(n, \mu) = (\Phi^{N,N+K}_{n,\lambda})^{-1} \Psi^{N,N+K}_{n,\lambda}.$$
\end{definition}

We observe that $I_{N,K}(n, \mu)$ takes the basis $\tilde w_t$ to the basis $v_t$ by construction. Moreover this isomorphism must by construction take the decomposition
$$ \hsu(n,\lambda) \cong \bigoplus_{\mu\in \Gamma^{n-1,\lambda}_{N,K}}\hsu(n-1,\mu) \otimes \hsu(\mu\otimes 1, \lambda)$$
to the decomposition
$$ \Vdag(\Sib^\infty_{n}, \mu_n, \lambda^\dagger)\cong \bigoplus_{\mu \in \Gamma_{N,K}^{n,\lambda}} \Vdag(\Sib^\infty_{n-1}, \mu_{n-1}, \mu^\dagger)\otimes \Vdag(\Sib^\infty_2, \mu,\Box, \lambda^\dagger)$$
such that $\tilde w_t \otimes \tilde w_{\mu,\lambda}$ is taken to $v_t \otimes v_{\mu,\lambda}$ by the very construction of the bases.

We note that we can inductively determine bases of $H^{SU(N)}_K(\mu\otimes n', \lambda)$ and of $ \Vdag(\Sib^\infty_{n+n'}, \mu_{n+n'}, \lambda^\dagger)$  by using the above bases and decompositions jointly with
$$ H^{SU(N)}_K(\mu\otimes n', \lambda) \cong  \bigoplus_{\nu \in \Gamma_{N,K}}  H^{SU(N)}_K(\mu\otimes (n'-1), \nu) \otimes H^{SU(N)}_K(\nu\otimes 1, \lambda).
$$

Since $I_{N,K}(n, \mu)$ takes preferred bases to preferred bases, we 
get the following theorem.

\begin{theorem}\label{decompI}
The decomposition
$$ \hsu( n+n', \lambda) \cong \bigoplus_{\mu\in \Gamma_{N,K}}\hsu( n,\mu) \otimes \hsu(\mu\otimes n', \lambda),$$
is by the isomorphism $I_{N,K}(n + n', \lambda)$ taken to the decomposition
$$ \Vdag(\Sib^\infty_{n+n'}, \mu_{n+n'}, \lambda^\dagger)\cong \bigoplus_{\mu \in \Gamma_{N,K}} \Vdag(\Sib^\infty_{n}, \mu_{n}, \mu^\dagger)\otimes \Vdag(\tilde\Sib^\infty_{n'}, \mu,\mu_{n'}, \lambda^\dagger).$$
Moreover preferred bases are taken to preferred bases by this isomorphism.
\end{theorem}

\subsection{The isomorphism for general box-labeled objects}
Since there for any object $\alpha$ in $H$ exist some $n$ and a braid $T_{\alpha,n}$ representing a morphisms in $H(\alpha,n)$, we get the following theorem. 
\begin{theorem}\label{genbox}
For all objects $\alpha$ in the category $H$, there is a unique isomorphism
$$
I_{N,K}(\alpha) : \hsu(\alpha,0 ) \ra \Vdag_{N,K}(\Sib_\alpha, \lambda_\alpha)
$$
such that
$$I_{N,K}(n) = I_{N,K}(n, 0)$$
and which makes the following diagram commutative 
$$
\begin{CD} \hsu(\alpha,0)@>I_{N,K}(\alpha)>> \Vdag_{N,K}(\Sib_\alpha, \lambda_\alpha)\\ 
@VV T_{\alpha,\beta} V @V\Vdag_{N,K}(\phi_{T_{\alpha,\beta}}) VV\\
 \hsu(\beta,0)@>I_{N,K}(\beta)>> \Vdag_{N,K}(\Sib_\beta, \lambda_\beta)
\end{CD}
$$
for all pairs of objects $\alpha,\beta$ in $H$ and all braid $T_{\alpha,\beta}$ representing morphisms in $H(\alpha,\beta)$.
\end{theorem}

From the above theorem, it is clear that  $\hsu(\alpha,\mu)$ and $\Vdag_{N,K}(\Sib_\alpha^\infty, \lambda_\alpha, \mu)$ are isomorphic irreducible $\hsu(\alpha,\alpha)$-module and that all of this algebra's irreducible modules are of this form by Theorem \ref{KanieTh} and Corollary \ref{cor}.

\subsection{The isomorphism for general labels}\label{igl}

Let  $\lambda_i \in \Gamma_{N,K}$. Let $n_i = |\lambda_i|$, $i= 1, \ldots, l$ and let $n = n_1+\ldots + n_l$. Let $\ell = (\ell_1, \ldots, \ell_l)$ be a framed set of points in the interior of $D^2$. Further let $\lambda$ be the object in $\hsu$ obtained by labelling $\ell_i$ by $\lambda_i$.

 Now let $\tilde \lambda = E(\lambda)$,
where $E(\lambda)$ is the object in $H$ defined in Section \ref{sec:C-compHecke-cat}.  We denote the labelling of all points in $\tilde \lambda$ simply also by $\tilde \lambda$. Further let $\tilde\lambda_i$ denote the object $E(\ell_i,\lambda_i)$ together with the labelling of all points in $\ell_i$ by boxes.

From Theorem \ref{genbox} we have the isomorphism
\begin{equation}\label{isoboxlambda}
I_{N,K}(\tilde \lambda) : \hsu(\tilde \lambda,0 ) \ra \Vdag_{N,K}(\Sib_{\tilde \lambda}, \tilde \lambda).
\end{equation}
By definition
$$
\hsu(\lambda,0) =   \pi_\lambda  \hsu(\tilde \lambda,0 ).
$$
By factorization in the boundaries of small disjoint discs $D_i$, $i=1, \ldots, l$ inside the unit disc, such that  $\ell_i \in D_i$, we get the isomorphism
$$
\Vdag_{N,K}(\Sib_{\ell}, \tilde \lambda) \cong \bigoplus_{\mu\in \Gamma_{N,K}^{\times l}} \Vdag_{N,K}(\Sib_\ell,\mu) \otimes \otimes_{i=1}^l \Vdag_{N,K}(\Sib_{\ell_i}^\infty,\tilde\lambda_i,\mu_i).
$$
We observe that this is a $\hsu(\tilde\lambda_1,\tilde\lambda_1)\times \ldots \times \hsu(\tilde\lambda_l,\tilde\lambda_l)$-module isomorphism.
By writing the identity in $\hsu(\tilde \lambda,\tilde \lambda)$ as a sum of minimal central idempotents in analogy with (\ref{1}) and inserting them just below $D_i\times 1$, $i=1,\ldots, l$, we also get the decomposition
\begin{equation}\label{fachom}
\hsu(\tilde \lambda,0) \cong \bigoplus_{\mu\in \Gamma_{N,K}^{\times l}} \hsu(\mu,0) \otimes \otimes_{i=1}^l \hsu(\tilde\lambda_i,\mu_i)
\end{equation}
which is also a $\hsu(\tilde\lambda_1,\tilde\lambda_1)\times \ldots \times \hsu(\tilde\lambda_l,\tilde\lambda_l)$-module isomorphism. Hence we see that 
the isomorphism (\ref{isoboxlambda}) must preserve these decompositions, thus we have proved the following Theorem.

\begin{theorem}\label{glhom}
For all objects $\lambda$ in $\hsu$ as above, there is a unique isomorphism
$$
I_{N,K}(\lambda) : \hsu(\lambda,0) \ra \Vdag_{N,K}(\Sib_\ell,\lambda) 
$$
and unique isomorphims 
$$
I_{N,K}(\tilde\lambda_i,\mu_i) : \hsu(\tilde \lambda,\mu_i) \ra \Vdag_{N,K}(\Sib_{\ell_i}^\infty,\tilde\lambda_i,\mu_i)
$$
for all $\mu_i\in \Gamma_{N,K}$ and $i=1,\ldots,l$ which for $\tilde \lambda_i=n_i$ agrees with the isomorphism from Definition \ref{hsuglabel}, are compatible with diffeomorphisms of the unit disc, that induces the identity on the boundary and such that under the above identifications
$$
I_{N,K}(\tilde \lambda) = \bigoplus_{\mu\in \Gamma_{N,K}^{\times l}} I_{N,K}(\mu) \otimes \otimes_{i=1}^l I_{N,K}(\tilde\lambda_i,\mu_i).
$$
\end{theorem}

\subsection{The isomorphism for arbitrary genus zero marked surfaces}\label{isog0}
Let 
$$\Sib = (\Si, P, V, L)$$ 
be a marked connected surface of genus zero and
$\lambda$ a labelling of it by labels from $\Gamma_{N,K}$, hence $\lambda : P \ra \Gamma_{N,K}$. Let $n=| P|$, then $\Sib_n$ is the standard surface of the same type as $\Sib$. Choose a morphism $\phi$ from $\Sib_n$ to $\Sib$.  This induces a labelling say $\lambda_0$ of $n$ which under $\phi|_{P_0}$ matches up with $\lambda$. Further let 
$$
\bar \lambda_0 = \otimes_{i=1} ^{|P_0|}\lambda_0(p_i).
$$
As was explained in section \ref{RTmf}, the modular functor $ \V_K^{SU(N)}$ assigns the vector space $ \V_K^{SU(N)}(\Sib,\lambda)$ represented by $\hsu(\bar \lambda_0,0)$ at the parametrization $(\Sib,\phi)$.

\begin{theorem}
There is a unique isomorphism
$$
I_{N,K}(\Sib, \lambda) : \V_K^{SU(N)}(\Sib,\lambda) \ra  \Vdag_{N,K}(\Sib,\lambda) 
$$
which for any parametrization $(\Sib,\phi)$ is represented by the isomorphism
$$
\Vdag_{N,K}(\phi) \circ I_{N,K}(\bar \lambda_0) : \hsu(\bar \lambda_0,0) \ra \Vdag_{N,K}(\Sib,\lambda).
$$
\end{theorem}
 
\begin{proof}
  We just need to show that $I_{N,K}(\Sib, \lambda) $ is well
  defined. Suppose $(\Sib,\phi_1)$ and $(\Sib,\phi_2)$ are two
  parametrizations. Then $\phi_1^{-1} \circ \phi_2$ is a morphism from
  $\Sib_n$ to it self, hence it can be represented by a braid
  $T_{\phi_1^{-1} \circ \phi_2}$ which represents an element in
  $\hsu(\bar \lambda_0,\bar \lambda_0)$. Now we have the commutative
  diagram
$$
\begin{CD} \hsu(\bar \lambda_0,0)@>I_{N,K}(\bar \lambda_0)>> \Vdag_{N,K}(\Sib_n, \lambda_0) @>\Vdag_{N,K}(\phi_1)>> \Vdag_{N,K}(\Sib, \lambda) \\
  @VV T_{\phi_1^{-1}\circ \phi_2} V @V\Vdag_{N,K}(\phi_1^{-1}\circ \phi_2) VV @V=VV\\
  \hsu(\bar \lambda_0,0)@>I_{N,K}(\bar \lambda_0)>>
  \Vdag_{N,K}(\Sib_n, \lambda_0) @>\Vdag_{N,K}(\phi_2)>>
  \Vdag_{N,K}(\Sib, \lambda)
\end{CD}
$$
which shows the induced isomorphism with respect to $\phi_1$ is the
same as the one induced from $\phi_2$.  
\end{proof}

We extend the isomorphism $I_{N,K}$ to disconnected surfaces of genus zero, by taking the tensor product of the isomorphisms for each component.

Suppose that $(\tilde P_p,\tilde V_p)$ is obtained from $\lambda(p)$ by the expansion $E$ introduced in section \ref{sec:C-compHecke-cat} for each $p\in P$. We let
$$
(\tilde P,\tilde V) = \bigsqcup_{p\in P} (\tilde P_p,\tilde V_p)
$$
Let $\tilde \lambda_p$ be the corresponding labeling of all points in $\tilde P_p$ by $\Box$ and similarly $\tilde \lambda$ assigns $\Box$ to all points in $\tilde P$. Let $\tilde \Sib = (\Si, \tilde P,\tilde V, L)$ and  $\tilde \Sib_p = (\Si, \tilde P_p\sqcup\{\infty\},\tilde V_p\sqcup \{v_\infty\}, L)$, where $v_\infty$ is the direction of the negative real axis at infinity.

\begin{theorem}\label{compgluebox}
We have the following commutative diagram 
$$
\begin{CD} \V_K^{SU(N)}(\tilde\Sib, \tilde \lambda) @>>> \bigoplus_{\mu: P \ra \Gamma_{N,K}} \V_K^{SU(N)}(\Sib,\mu) \otimes \otimes_{p\in P} \V_K^{SU(N)}(\Sib_p,\tilde\lambda_p,\mu(p))\\ 
@V I_{N,K}(\tilde\Sib, \tilde \lambda) VV @V \bigoplus_{\mu: P \ra \Gamma_{N,K}} I_{N,K}(\Sib,\mu) \otimes \otimes_{p\in P} I_{N,K}(\Sib_p,\tilde\lambda_p,\mu(p))VV \\
 \Vdag_{N,K}(\tilde\Sib, \tilde \lambda) @>>> \bigoplus_{\mu: P \ra \Gamma_{N,K}} \Vdag_{N,K}(\Sib,\mu) \otimes \otimes_{p\in P} \Vdag_{N,K}(\Sib_p,\tilde\lambda_p,\mu(p))
\end{CD}
$$
where the horizontal arrows are the factorization isomorphisms.
\end{theorem}

\begin{proof}
  Since the factorization isomorphisms is compatible with
  morphisms of labelled marked surface, we can pick an morphism from
  $\tilde \Sigma$ to the standard surface of the same type, which is
  $\Sib_n$ for some $n$ and then simply just check the corresponding
  diagram for $\Sib_n$ factored along a number of disjoint embedded
  regular discs with centres on the real line of the appropriate
  radii. For for this standard surface and these standard
  factorization curves, the commutativity of the diagram in this
  theorem follows immediately by the remark just above equation
  (\ref{facc}) and repeated use of Theorem \ref{decompI} and Lemma
  \ref{obmf}.  
\end{proof}

\subsection{Compatibility with glueing in genus zero}

Let $\Sib = \Sib^+ \sqcup \Sib^-$, where 
$$\Sib^\pm = (\Si^\pm, \{p^\pm\}\sqcup P^\pm,\{v^\pm\}\sqcup V^\pm,L)$$ 
are marked surfaces of genus zero. Let $\lambda^\pm : P^\pm \ra \Gamma_{N,K}$ be a labeling of $P^\pm$. 
Further let $\Sib_c$ be the glueing of $\Sib$ with respect to $(p^\pm,v^\pm)$ and $\lambda_c$ the corresponding labeling of $P = P^-\sqcup P^+$.
\begin{theorem}\label{glueAx}
We have the following commutative diagram
$$
\begin{CD} \V_K^{SU(N)}(\Sib_c,\lambda_c)  @>>>  \bigoplus_{\mu\in \Gamma_{N,K}}\V_K^{SU(N)}(\Sib^-,\mu, \lambda^-)\otimes \V_K^{SU(N)}(\Sib^+,\mu^\dagger, \lambda^+)\\ 
@V I_{N,K}(\Sib_c,\lambda_c) VV @V \bigoplus_{\mu\in \Gamma_{N,K}}I_{N,K}(\Sib^-,\mu, \lambda^-) \otimes  I_{N,K}(\Sib^+,\mu^\dagger, \lambda^+)VV \\
\Vdag_{N,K}(\Sib_c,\lambda_c)  @>>> \bigoplus_{\mu\in \Gamma_{N,K}}\Vdag_{N,K}(\Sib^-,\mu, \lambda^-)\otimes \V_{N,K}(\Sib^+,\mu^\dagger, \lambda^+)
\end{CD}
$$
where the horizontal arrows are the factorization isomorphisms.
\end{theorem}

\begin{proof}
  Let $\tilde\Sigma^\pm$ be the expansion of $\Sigma^\pm$
  determined by $E(\lambda^\pm)$ and let $\tilde \Sigma_c$ be the
  resulting expansion of $\Sigma_c$.  Again the invariance of
  factorization under morphisms of surfaces allows us to assume that
  $\tilde\Sigma^\pm$ are standard. Further we can assume that $\tilde
  \Sigma_c$ also has been identified with the standard surface of its
  type. Then the commutativity follows by the remark just above
  quantum (\ref{facc}) and repeated use of Theorem \ref{decompI} and
  Lemma \ref{obmf}, as in the previous proof.  
\end{proof}

\subsection{The isomorphism of the genus zero part of the modular functors}
We summarise our construction so fare in the following Theorem.

\begin{theorem}\label{0=g}
We have an isomorphism $I_{N,K}$ from the genus zero part of the modular functor $\V_K^{SU(N)}$ to the genus zero part of the modular functor $\Vdag_{N,K}$.
\end{theorem}

This theorem follows directly from the fact that this isomorphisms per construction is compatible with morphisms of marked surfaces, the definition of the isomorphisms for disconnected surfaces and then from Theorem \ref{glueAx}.

\section{The $S$-matrices and the higher genus isomorphism}\label{SMatrixhg}

First we recall the main theorem from \cite{AU3}. 
\begin{theorem}\label{0=>g}
Suppose $\V_i$, $i=1,2$ are modular functor and we have isomorphisms
$$
I(\Sib,\lambda) : \V_1(\Sib,\lambda) \ra \V_2(\Sib,\lambda)
$$
for all genus zero labeled marked surfaces $(\Sib,\lambda)$, which is compatible with disjoint union and glueing within genus zero labeled marked surfaces, then there exists a unique extension of $I$ to all labeled marked surfaces of all genus, which gives a full isomorphism 
$$
I : \V_1 \ra \V_2
$$
of modular functors.
\end{theorem}

We recall that this Theorem is proved by showing that the two $S$-matrices of the two theories $\V_i$, $i=1,2$ agree. Once we have this, it is clear that there is a unique isomorphism $I(\Sib,\lambda)$, which is compatible with obtaining $\Sib$ as the glueing of trinions, up to morphism of labeled marked surfaces, as explained in \cite{AU3}.

Our main Theorem \ref{Main} follows now directly from Theorem \ref{0=g} and Theorem \ref{0=>g}. We remark that Theorem \ref{0=>g} does not require the modular functors in question to have duality. In fact at present we do not have a geometric construction of a duality structure for $ \Vdag_{N,K}$. That $ \Vdag_{N,K}$ has one is a consequence of our main Theorem \ref{Main}.

\end{document}